\documentclass[%
aip,
superscriptaddress,
amsmath,amssymb,
reprint,
]{revtex4-2}

\usepackage{graphicx}
\usepackage{dcolumn}
\usepackage{bm}
\usepackage{theoremref}
\usepackage{mathtools}
\newcommand\dataset[1]{\textsc{\texttt{#1}}}
\newtheorem{theorem}{Theorem}
\numberwithin{theorem}{section}

\newtheorem{problem}[theorem]{Problem}

\usepackage{algorithm}
\usepackage[noend]{algpseudocode}
\usepackage{float}
\usepackage[utf8]{inputenc}
\usepackage[T1]{fontenc}
\usepackage{booktabs}
\usepackage{amsmath}
\usepackage{mathptmx}

\usepackage{graphicx}
\usepackage{dcolumn}
\usepackage{xcolor}
\usepackage{bm}
\usepackage{array}

\makeatletter
\def\namedlabel#1#2{\begingroup
	\def\@currentlabel{#2}%
	\label{#1}\endgroup
}
\makeatother

\def\1{\mathbf{1}}

\def\R{\mathbb{R}}

\def\P{\mathbb{P}}
\def\F{\mathcal{F}}
\def\eps{\varepsilon}

\def\E{\mathbb{E}}

\def\x{\mathbf{x}}
\def\X{\mathcal{X}}

\def\1{\mathbf{1}}

\def\Beta{\boldsymbol{\beta}}

\def\W{\mathbf{W}}
\def\w{\mathbf{w}}

\def\x{\mathbf{x}}
\def\h{\mathbf{h}}

\def\Beta{\boldsymbol{\beta}}

\DeclareMathOperator{\MAT}{\texttt{\textsc{MAT}}}
\DeclareMathOperator{\VEC}{\texttt{\textsc{VEC}}}
\DeclareMathOperator*{\argmax}{arg\,max}
\DeclareMathOperator*{\argmin}{arg\,min}

\usepackage{etoolbox}
\makeatletter
\patchcmd{\@thm}{\thm@headfont{\scshape}}{\thm@headfont{\scshape\bfseries}}{}{}
\patchcmd{\@thm}{\thm@notefont{\fontseries\mddefault\upshape}}{}{}{}
\makeatother

\makeatletter
\renewenvironment{description}
{\list{}{\leftmargin=5pt
		\labelwidth\z@ \itemindent-\leftmargin
		}}%
{\endlist}
\makeatother

\usepackage{hyperref}
\hypersetup{colorlinks,linkcolor={red},citecolor={blue},urlcolor={red}}

\begin{document}
	
	\preprint{AIP/123-QED}
	
	\title{A latent linear model for nonlinear coupled oscillators on graphs}
	${}$
	\vspace{0.5cm}
	
	\author{Agam Goyal}
	\thanks{Co-first author}
	\affiliation{
		Department of Computer Science, University of Wisconsin - Madison, WI 53706}
	\affiliation{
		Department of Mathematics, University of Wisconsin - Madison, WI 53706}
	
	\author{Zhaoxing Wu}
	\thanks{Co-first author}
	\affiliation{
		Department of Statistics, University of Washington, Seattle, WA 98195}
	
	\author{Richard P. Yim}\affiliation{
		Bakar Computational Health Sciences Institute, University of California, San Francisco, CA 94143 \looseness=-1}
	
	\author{Binhao Chen}\affiliation{
		Department of Computer Science, Brown University, Providence, RI 02912}
	
	\author{Zihong Xu}\affiliation{
		Department of Mathematics, University of Wisconsin - Madison, WI 53706}
	
	\author{Hanbaek Lyu}
	\thanks{Corresponding author; \texttt{hlyu@math.wisc.edu}}
	\affiliation{
		Department of Mathematics, University of Wisconsin - Madison, WI 53706}
	
	\begin{abstract}
		A system of coupled oscillators on an arbitrary graph is locally driven by the tendency to mutual synchronization between nearby oscillators, but can and often exhibit nonlinear behavior on the whole graph. Understanding such nonlinear behavior has been a key challenge in predicting whether all oscillators in such a system will eventually synchronize. In this paper, we demonstrate that, surprisingly, such nonlinear behavior of coupled oscillators can be effectively linearized in certain latent dynamic spaces. The key insight is that there is a small number of `latent dynamics filters', each with a specific association with synchronizing and non-synchronizing dynamics on subgraphs so that any observed dynamics on subgraphs can be approximated by a suitable linear combination of such elementary dynamic patterns. Taking an ensemble of subgraph-level predictions provides an interpretable predictor for whether the system on the whole graph reaches global synchronization. We propose algorithms based on supervised matrix factorization to learn such latent dynamics filters. We demonstrate that our method performs competitively in synchronization prediction tasks against baselines and black-box classification algorithms, despite its simple and interpretable architecture. 
	\end{abstract}
	
	\maketitle
	
	\section{Introduction}\label{sec:introduction}
	
	If a group of people is given local clocks with arbitrarily set times, and there is no global reference (for example, GPS), is it possible for the group to synchronize all clocks by only communicating with nearby members? In order for a distributed system to be able to perform high-level tasks that may go beyond the capability of an individual agent, the system must first solve a ``clock synchronization'' problem to establish a shared notion of time. The study of synchronization of coupled oscillators has been an important subject of research in mathematics and various areas of science for decades\cite{strogatz2000kuramoto, acebron2005kuramoto}, with fruitful applications in many areas, including wireless sensor networks, wildfire monitoring, electric power networks, robotic vehicle networks, and large-scale information fusion \cite{nair2007stable, pagliari2011scalable, dorfler2012synchronization}. 
	
	A system of coupled oscillators is said to be (globally) \textit{synchronized} if all oscillators are at a consensus in terms of their phase or oscillation frequency. In this work, we consider oscillators of identical frequencies and only phase synchronization. Such a global state may or may not emerge depending on how the oscillators interact along the edges of the graph, how such local interaction leads to larger-scale interactions, and so on. In spite of several sufficient conditions on model parameters (e.g., large coupling strength \cite{dorfler2011critical}) or initial configuration (e.g., phase concentration \cite{klinglmayr2012guaranteeing} within an open semicircle), it is usually analytically intractable to predict whether a given system of coupled oscillators with arbitrary underlying graph structures will eventually synchronize, more so when the underlying graph is heterogeneous and the initial phase configuration is not confined in a small arc of the phase space. Furthermore, the interplay between the nonlinear dynamics and network topology can often give rise to highly nonlinear phenomena\cite{skardal2020higher, d2019explosive}, which makes it intriguingly hard to study and understand their properties.
	
	With the revolutionary success of machine learning methods in various tasks such as image classification and natural language processing, there has been a surge of interest in employing these methods to study scientific problems that have been previously believed to be extremely difficult \cite{jumper2021highly, koren2009matrix, mehta2021alx}. This is also the case for the problem of synchronization prediction, where a number of tools in machine learning have been applied to study the properties of coupled oscillator systems \cite{fan2021anticipating, guth2019machine, chowdhury2021extreme}. The recent work of Bassi et al. \cite{bassi2022learning} in particular demonstrated that after a proper reformulation, the synchronization prediction problem on randomly generated graphs can be effectively solved by training binary classification algorithms on a large dataset of synchronizing and non-synchronizing examples, in the sense that the resulting prediction accuracy significantly outperforms a baseline predictor that uses the concentration principle in coupled oscillator theory \cite{nishimura2011robust, klinglmayr2012guaranteeing, lyu2018global}. 
	
	However, a key question that remains unanswered is whether we can actually gain any scientific insight into coupled oscillator systems from the advantages of analyzing a massive amount of data using machine learning methods. For instance, we would like to understand `how' a well-performing model is able to make these predictions on coupled oscillator systems and what features in graphs or dynamics it considers most important for this task. Toward this goal, in this paper, we propose an \textit{interpretable} model for synchronization prediction that we call the \textit{Latent Linear Dynamics Model} (LLDM). The logic behind the model is very simple: 
	\textit{Given an observed dynamics on a subgraph, first compute `proximity scores' for how much a set of prescribed patterns we observe there, and then use the proximity scores for those patterns as an input to a logistic classifier.}

	The key challenge in our approach is to figure out what fundamental patterns of dynamics on subgraphs we seek to observe for the purpose of synchronization prediction. 
	For example, if the underlying graph is very dense, then it will be likely that the dynamics will eventually synchronize. Also, if the observed phase configuration is confined in an open semicircle, then we know the system will eventually synchronize. If the graph is sparse and contains a long cycle, then it would be hard to see eventual synchronization. While such `patterned behaviors' are informed by the existing knowledge on coupled oscillators, our novel approach here is to learn such `critical patterns for synchronization' directly from the data. \\
	
	We summarize our key contributions through this work:
	\begin{enumerate}
		\item[1.] We propose a novel and interpretable framework for the prediction of synchronization in coupled oscillators, leveraging the feature learning capabilities of matrix factorization techniques to learn latent linear representations of underlying network dynamics.
		\item[2.] We propose various ways of approaching the synchronization prediction problem using \textup{LLDM}, by the use of data-informed and computationally efficient, theory-informed approaches, in addition to using both supervised and unsupervised matrix factorization techniques.
		\item[3.] We propose a compute-efficient and novel method for the prediction of dynamics synchronization on large-scale graphs by the use of \textup{LLDM} on a set of subgraphs sampled by motif sampling techniques, followed by recursive averaging of predicted probabilities.
	\end{enumerate}
	
	To the best of our knowledge, this is the first work to study the synchronization of small and large-scale coupled oscillator dynamical systems on graphs through the lens of representation learning techniques that also focus on interpretability--- a critical aspect of modeling dynamical systems.
	
	\subsection{Related Works}\label{sec:related_works}
	
	\subsubsection{Machine Learning for Synchronization Prediction}
	
	There has been a surge in studies that have approached the study the dynamic oscillator systems by using machine learning techniques. Some studies such as Thiem et al.\cite{thiem2020emergent} used Feed Forward Neural Networks\cite{bishop2006pattern} (FFNNs) to analyze Kuramoto dynamics in specific, while Hefny et al.\cite{hefny2015supervised} use LASSO regression for modeling independent subsystems of dynamical systems. Itabashi et al.\cite{itabashi2021evaluating} use features derived from early-stage topological dynamics to classify Kuramoto oscillator dynamics. Furthermore, Bassi et al.\cite{bassi2022learning} show that various classical machine learning algorithms can be used for the synchronization prediction problem by training them on a large dataset of synchronizing and non-synchronizing coupled oscillator systems on randomly generated graphs. Their method was applied to the Kuramoto oscillators as well as discrete oscillator models such as the Firefly Cellular Automata (FCA\cite{lyu2015synchronization}) and Greenberg-Hastings Model (GHM\cite{greenberg1978spatial}). Recently, Chen et al.\cite{chen2023optimal} proposed to use reinforcement learning to find an optimal pulse-interaction mechanism that optimizes the probability of synchronization of pulse-coupled oscillators, while Mahlow et al.\cite{mahlow2023predicting} proposed to utilize $k$-nearest neighbor regressor to predict the emergence of environment-induced spontaneous quantum synchronization in an open system setting.
	
	Despite the increasing interest in utilizing these methods to study dynamical systems, there remains a gap in understanding what features are crucial for these models to make predictions. We aim to bridge this gap in our work, where we leverage feature-representation learning techniques like non-negative matrix factorization\cite{lee1999learning} and supervised matrix factorization\cite{mairal2008supervised, lee2023supervised} (See Appendices \ref{app: nmf} and \ref{app: SMF} for details) to provide an interpretable framework for synchronization prediction.  
	
	\subsubsection{Matrix Factorization Techniques and interpretable feature extraction}
	
	Matrix factorization techniques have proven to be powerful tools for describing various latent features of data of interest in terms of a `linear combination' of atomic elements. This problem has been studied for many decades and has been used in various scientific fields\cite{berry2005email, boutchko2015clustering, sitek2002correction, berry2005email, berry2007algorithms, chen2011phoenix, taslaman2012framework, boutchko2015clustering, ren2018non}. More generally, low-dimensional feature extraction techniques have been used extensively in the last few decades for complex tasks that involve studying the local interactions of elements in various fields. Highly accurate and precise reconstruction of billions of protein structures \cite{jumper2021highly}; OTT media recommendation systems decoding users' item-response patterns \cite{koren2009matrix}; the innovation in novel TPU architectures to allow efficient matrix factorization \cite{mehta2021alx}; and dating back to the PageRank \cite{brin1998anatomy} search algorithm for ranking internet websites and web pages; these are all some of the most well-known and important applications of such techniques. 
	
	Despite enjoying fruitful applications in the aforementioned areas, the potential of these techniques has not been harnessed widely in the context of coupled oscillators. Recently, Luo\cite{luo2023unraveling} proposed to decompose network dynamics into a composite of weighted principal components, and subsequently learn the governing differential equations using sparse regression.
	In this work, we make use of matrix factorization-based approaches in conjunction with a Markov-chain Monte-Carlo subgraph sampling algorithm \cite{lyu2023sampling} to learn underlying features from the data and make predictions regarding the synchronization of oscillator dynamic systems.

	\begin{figure*}[ht]
		\centering
		\includegraphics[width=1\textwidth]{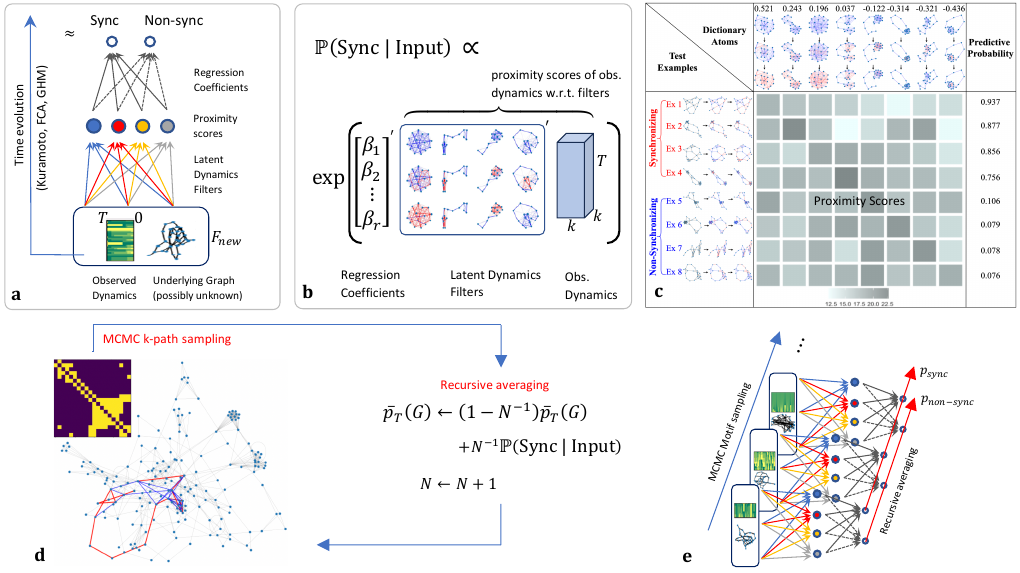}
		\caption{(Panels \textbf{a}-\textbf{c}) Scheme of the Latent Dynamics Model (LLDM) for synchronization prediction. Input dynamics on a $k$-node graph observed for $T$ iterations is represented as a $k\times k  \times T$ tensor.  Taking convolution of the input dynamics with $R$ (4 in the figure), latent dynamics filters (nonnegative tensors of shape $k \times k\times T$) give proximity scores for $R$ patterns, which are combined into one scalar for a final score for the predictive probability of eventual synchronization. The predictive probability of synchronization is proportional to the exponential of the final score. (Panels \textbf{d}-\textbf{e}) A uniformly randomly sampled $20$-path (with red edges) and additional edges on the induced subgraph (in blue) with the corresponding adjacency matrix on the top left. We recursively compute the predictive probabilities of synchronization using dynamics on such sampled subgraphs, and the average value converges to the predictive probability of synchronization of the whole graph, which is the expectation of the predictive probabilities over subgraphs induced on $k$-paths.  } 
		\label{fig:GLLM_scheme}
	\end{figure*}

	\section{Statement of the problem}\label{sec:problem_statement}
	
	A graph $G=(V,E)$ consists of sets $V$ and $E$, node set and edge set, respectively. Let $\Omega$ denote the \textit{phase space} of each node, which may be taken to be the unit circle $\mathbb{R}/2\pi \mathbb{Z}$ for continuous-state oscillators and the discrete circle $\mathbb{Z}/\kappa \mathbb{Z}$, $\kappa\in \mathbb{N}$ for finite-state oscillators. We call a map $X:V\rightarrow \Omega$ a \textit{phase configuration}, and say it is \textit{synchronized} if it takes a constant value across nodes (i.e., $X(v)=Const.$ for all $v\in V$). A \textit{coupling} is a function $\mathcal{C}$ that maps each pair $(G,X_{0})$ of graph and initial configuration $X_{0}:V\rightarrow \Omega$ deterministically to a \textit{trajectory} $(X_{t})_{t\ge 0}$ of phase configurations $X_{t}:V\rightarrow \Omega$. In this paper,  we consider $\mathcal{C}$ to be the time evolution rule for Kuramoto model\cite{kuramoto2003chemical}, Firefly Cellular Automoata (FCA)\cite{lyu2015synchronization, lyu2023time}, and Greenberg-Hastings Model\cite{greenberg1978spatial}. We use a discretization of Kuramoto model (see \eqref{eq:kura_discretize} in the appendix) and an `iteration' of Kuramoto dynamics refers to applying one step of the difference equation. 
	See Appendix \ref{app: oscillators} for details on each of these coupled oscillator models.
	
	The main problem we investigate in this work is to predict the synchronization of coupled oscillators on a large-scale graph $G$ using subgraph-level information. This means that we observe some subgraphs of fixed size $k$  (potentially much smaller than the number of nodes in $G$) and the dynamics in $G$ are restricted on these subgraphs. In order to make this setting precise, we formulate the following sampling oracle: 
	\begin{description}
		\item[Subgraph Sampling Oracle.] \textit{Given a graph $G$ and a fixed integer $k\ge 0$, we can sample a $k$-node connected subgraph $H$ of $G$ and observe dynamics on $G$ restricted on $H$ up to a fixed number of iterations $T_{0}$. However, we cannot observe dynamics on larger subgraphs and more than $T_{0}$ iterations.}
	\end{description}
	
	We can now formulate the main problem we aim to address in this work: \textit{Large-scale synchronization prediction by local dynamics decomposition}:
	\begin{problem} 
		\label{problem:main2}
		Let $(X_{t})_{t\ge 0}$ be a coupled oscillator dynamics on a (possibly large) connected graph $G$ governed by a coupling $\mathcal{C}$. Suppose that we have a sampling oracle for subgraph size $k$ and time horizon $T_{0}$. 
		\begin{description}
			\item[(i)] (Dynamics decomposition) Observed dynamics on subgraphs can be approximately decomposed into a linear combination of some key dynamics patterns. 
			
			\item[(ii)] (Synchronization prediction) Using the decomposition in \textbf{(i)}, one can predict the following indicator variable: 
			\begin{align}
				\mathbf{1}(\text{$X_{t}$ on $G$ is eventually  synchronized}).
			\end{align}

		\end{description}

	\end{problem}
	
	Our goal is to use SMF to learn low-rank latent factors that offer interpretable, data-reconstructive, and class-discriminative features, addressing the challenges posed by high-dimensional data.
	
	Our approach has two components. First, we will develop an interpretable model to predict the synchronization indicator at the $k$-node subgraph level. Second, we apply the trained subgraph-level model to a randomly sampled $k$-node along with the observed dynamics. The expectation of the predicted probability over random $k$-node subgraphs will be the predictive probability for the parent graph to eventually synchronize.

	\section{Model description}
	
	At a high level, our model predicts the eventual synchronization on the whole graph $G$ by averaging the predictions on many suitably chosen subgraphs of $F$. Below, we first describe how we make synchronization prediction using a single subgraph and then discuss how to efficiently combine subgraph-level prediction with our particular choice of the subgraph sampling oracle. 
	
	\subsection{Latent linear dynamics model on a single subgraph}
	\label{sec:define_model}
	
	Suppose that we have a system of coupled oscillators on a connected graph $G=(V,E)$. The goal of LLDM is to model the predictive probability that the system will eventually synchronize based on the observation of dynamics up to $T$ iterations restricted on a $k$-node subgraph of $G$, say $F$. (Here we assume $k\le |V|$ and allow $F=G$.) Since observing dynamics on a subgraph during a fixed time period only gives partial information on the long-term dynamics on the whole graph $G$, we model the indicator variable that the dynamics on $G$ will eventually synchronize given this partial information as a Bernoulli random variable with unknown success probability. 
	
	For a precise formulation, let $(X_{t})_{0\le t < T}$, $X_{t}:V\rightarrow \mathbb{Z}/\kappa\mathbb{Z}$ denote the dynamics on $G$ that is assumed to have evolved according to some coupled oscillator model (e.g., FCA, Kuramoto, or GHM). Let $X_{t}[F]$ denote the restriction of $X_{t}$ on the node set of $F$. Our basic modeling assumption is the following: 
	\begin{align}
		&\mathbf{1}\left( \textup{$X_{t}$ synchronizes as $t\rightarrow\infty$}\,|\,  F, (X_{t}[F])_{0\le t < T}   \right) \nonumber \\ 
		&\qquad \sim \text{Bernoulli}(p_{T}(F)),
	\end{align}
	where $p_{T}(F)$ is the unknown probability of eventual synchronization of the dynamics on $G$ given the partial observation.

	\begin{figure}[!ht]
		\centering
		\includegraphics[width=0.48\textwidth]{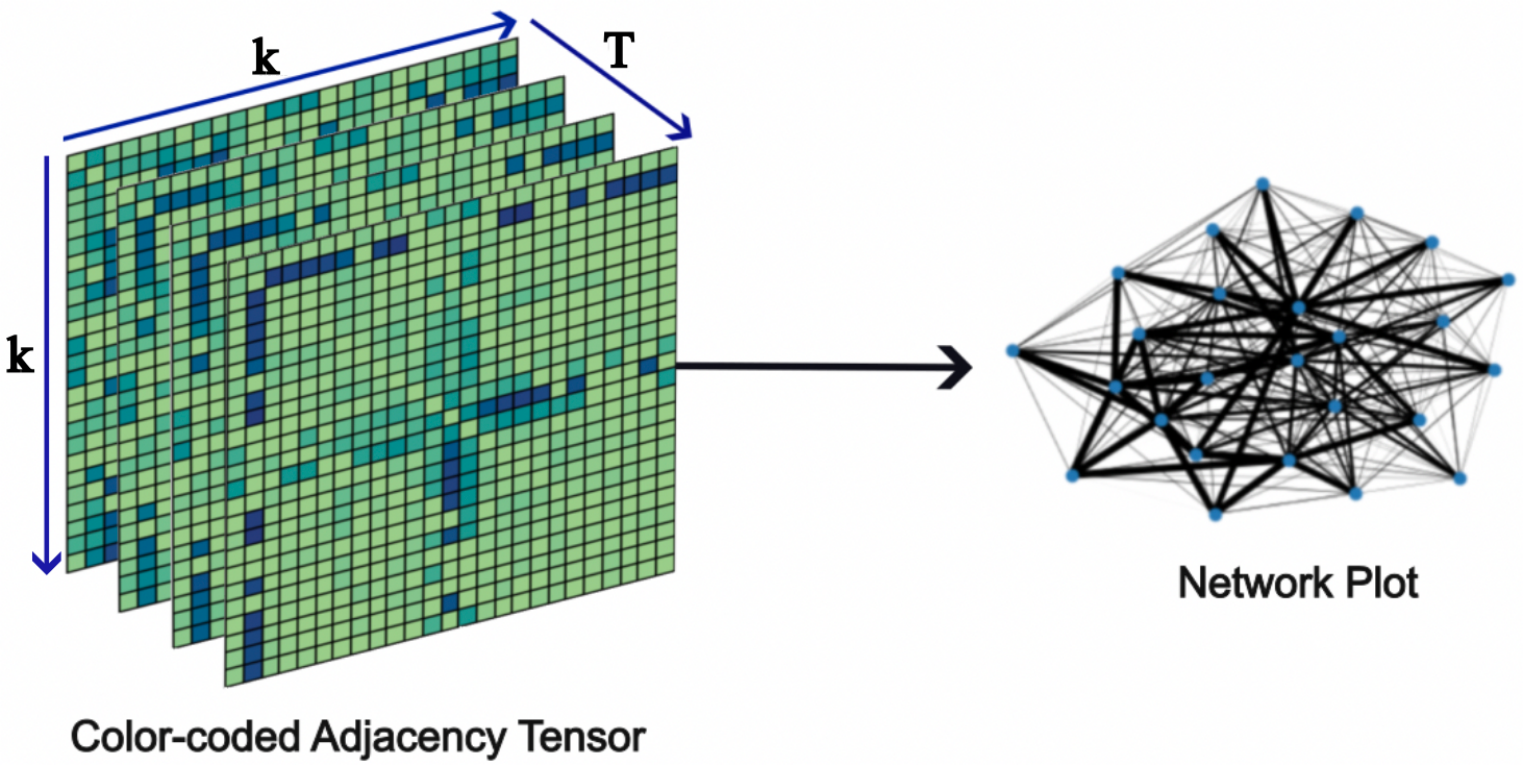}
		\caption{Example of a $k \times k \times T$ colored adjacency tensor (CAT) and the top slice as a graph network plot} 
		\label{fig:CAT}
	\end{figure} 
	
	Next, we introduce the key modeling assumption for LLDM that $p_{T}(F)$ above depends on a certain $R$-dimensional \textit{proximity score vector} $\mathbf{h}$ that is computed as follows. First, we represent the pair $(F, (X_{t}[F])_{0\le t < T})$ of input data by a nonnegative tensor $\X$ of shape $k\times k\times T$, where each slice $\X[:,:,t]$ represents the graph topology $F$ decorated by the phase configuration $X_{t}$ defined as 
	\begin{align}\label{eq:CAT}
		\X[i,j,t] := A_{ij} \min \{ \left( X_{t}(i) - X_{t}(j) \textup{(mod $\kappa$)},\right. \nonumber \\ 
		\left.X_{t}(j) - X_{t}(i) \textup{(mod $\kappa$)} \right) \},.
	\end{align}
	Here, $A=(A_{ij})$ is the adjacency matrix of $F$ defined by $A_{ij}:=\mathbf{1}(\{i,j\})\in E(F))$. We call the tensor $\X$ the \textit{colored adjacency tensor} (CAT) for $(F, (X_{t}[F])_{1\le t < T})$ (see Figure \ref{fig:CAT}). In Appendix \ref{app:graph_embedding}, we demonstrate that the use of various encoding methods of graph topology of $F$, as simple as the adjacency matrix or as advanced as modern graph embedding methods such as \textup{DeepWalk \cite{perozzi2014deepwalk}, graph2vec \cite{narayanan2017graph2vec}}, provides little to no improvement in prediction accuracy for the synchronization problem at hand. 
	
	Second, a key notion we introduce in this work is the `dictionary' of \textit{latent dynamics filters}, which is an $R-$tuple  $\mathcal{D}=[\F_{1},\dots,\F_{R}]$ of nonnegative tensors $\F_{i}\in \R_{\ge 0}^{k\times k\times T}$ of unit Frobenius norm. 
	Here $R$ is an integer parameter called the `size' of the dictionary $\mathcal{D}$. Each filter $\F_{i}$ represents an elementary pattern of coupled oscillator dynamics on a general $k$-node graph for $T$ iterations. We denote by $\langle \cdot, \cdot\rangle$ the inner product between two tensors of the same shape, which is the sum of the products of all corresponding entries. Since both the input tensor $\X$ and the filter $\F_{i}$ are assumed to be nonnegative, their inner product $\langle \F_{i}, \X \rangle$ can be interpreted as the \textit{proximity score} of the input dynamics $\X$ with respect to the dynamics pattern represented by the filter $\F_{i}$. Since there are $R$ filters, we can compress a single tensor data element, $\X$, into an $R$-dimensional vector of proximity scores 
	\begin{align}
		\mathbf{h}&:=\MAT(\mathcal D)^{T}\VEC(\mathcal X)  \nonumber \\
		&=  (\langle \F_{1}, \X \rangle,\dots,\langle  \F_{R}, \X \rangle)^{T},
	\end{align}
	where $\MAT(\cdot)$ and $\VEC(\cdot)$ are the matricization and the vectorization operators---fixing the lexicographic ordering of entries.  (Note that $\MAT(\mathcal{D})\in \R^{k^{2}T \times r}$ and $\VEC(\X)\in \R^{k^{2}T}$ so that $\MAT(\mathcal{D})^{T} \VEC(\X) \in \R^{r}$.) The $i$th coordinate of $\mathbf{h}$ measures how similar the observed dynamics on the subgraph $F$, encoded in the tensor $\mathcal{X}$, to the particular dynamics and networks encoded in the $i$th latent dynamic filter $\mathcal{F}_{i}$. 
	
	Now we suppose the probability $p_{T}(F)$ of eventual synchronization given data $\mathcal{X}\in \R^{k\times k\times T}$ is modeled as follows: 
	\begin{align}\label{eq:predictive_prob_LLDM}
		p_{T}(F) 
		= \frac{\exp\left(\mathbf \Beta^{T}\h \right)}{1
			+\exp\left(\Beta^{T} \h \right)} \,\,\,\, \textup{or} \,\,\,\, \textup{logit}(p_{T}(F)) = \Beta^{T} \h, 
	\end{align} 
	where $\Beta=(\beta_{1},\dots,\beta_{R})^{T}\in \R^{r}$ is a vector of regression coefficients and $\textup{logit}(p):=\frac{p}{1-p}$. Thus, LLDM for fixed subgraph size $k$ is parameterized by $(\mathcal{D},\Beta)$. From the above representation, we see that LLDM is a linear model on the latent space of features measured by the proximity score matrix. 
	
	Our general scheme of predicting synchronization of coupled oscillators using \textup{LLDM} is depicted in Figure \ref{fig:GLLM_scheme}. Figure \ref{fig:GLLM_scheme}\textbf{c} depicts eight observed FCA dynamics on subgraphs of \dataset{NWS} (see Table \ref{tab:networks}) in the rows (``Test Examples'') and eight latent dynamics filters with their corresponding regression coefficients in the columns (``Dictionary Atoms''). Each of them is $k\times k\times T$ tensors for $k=20$ and $T=50$, and three snapshots at times $0, 25, 50$ are shown with arrows indicating time evolution. The proximity scores are shown in the heat map (in grayscale) and the corresponding predictive probabilities (see \eqref{eq:predictive_prob_LLDM}) are shown in the last column. ``Ex1'' in Figure \ref{fig:GLLM_scheme}\textbf{c} is a synchronizing test example, which has the largest proximity score with the first latent dynamic filter (with regression coefficient 0.521), and has a large predictive probability of 0.917 for eventual synchronization. Also, ``Ex7'' is a non-synchronizing test example that has large proximity scores with the fifth ($-0.122$) and the seventh ($-0.321$)  latent dynamic filter and has a small predictive probability of 0.078 for eventual non-synchronization.

	\subsection{Choosing the sampling oracle: $k$-path motif sampling}\label{sec:path_sampling}
	
	In the previous section, we introduced LLDM with parameters $\mathcal{D}$ and $\Beta$, without any assumption on the subgraph $F$ on which we observe the dynamics on $G$. In order for this model to be effective, especially when $G$ is large and sparse, we may need to restrict the class of `appropriate' $k$-node subgraphs in $G$ sampled by our sampling oracle.  For instance, if we consider all induced subgraphs obtained by sampling $k$ nodes uniformly at random from $G$, then when $G$ is sparse, most of such subgraphs will be disconnected and have a few edges, so the dynamics observed on such subgraphs will not be informative of the dynamics on the whole graph $G$. Furthermore, LLDM assumes that the observed subgraphs come with prescribed node ordering so that their adjacency matrix, and in turn their CAT. Thus it is computationally beneficial to restrict ourselves to consider $k$-node connected subgraphs that have canonical node ordering. 
	
	We propose to consider $k$-node connected subgraphs that are obtained by first uniformly randomly sampling a `$k$-path' in the graph $G$ and then taking the induced subgraph on the sampled paths (i.e., including all edges between the sampled nodes in $G$). See Figure~\ref{fig:GLLM_scheme}\textbf{d} for an illustration. Here a \textit{$k$-walk} is a sequence $\x = (x_{1}, \ldots,x_{k})$ of $k$ nodes (which may or may not be distinct), such that $x_{j}$ and $x_{j+1}$ are adjacent for all $j \in \{1, \ldots, k-1\}$. A $k$-walk is a \textit{$k$-path} if all nodes in the walk are distinct. This sampling method has two notable advantages. First, it guarantees that the sampled $k$-node induced subgraph is connected with the minimum number of imposed edges. Second, it induces a canonical node ordering of the sampled subgraphs. In order for efficient sampling of a large number of $k$-paths approximately uniformly at random, we use the $k$-walk motif-sampling algorithm in Lyu et al.\cite{lyu2023sampling} (which is a Markov chain Monte Carlo (MCMC) algorithm) in conjunction with rejection sampling. This subgraph sampling oracle has been recently used in Lyu et al. \cite{lyu2021learning} for mesoscale network reconstruction.

	\subsection{Averaging subgraph-level predictions over many subgraphs} \label{sec:subgraph_to_parent}

	Using the $k$-path motif sampling method introduced in the previous section, we propose the following simple procedure to improve LLDM by averaging over many subgraphs. Namely, instead of using a single $k$-node subgraph $F$ and the dynamics observed on it to predict the synchronization indicator on $G$, we use many such subgraphs $F$ and average the corresponding predictive probabilities. This effectively combines subgraph-level predictions over many subgraphs. Accordingly, we define 
	\begin{align}\label{eq:global_predictive_prob}
		\overline{p}_{T}(G) :=  \E_{F}\left[ p_{T}(F) \right], 
	\end{align}
	where $F$ is a random $k$-node connected subgraph in $G$ induced on a uniformly random $k$-path in $G$, and $p_{T}(H)$ is the predictive probability (using LLDM, see \eqref{eq:predictive_prob_LLDM}) of $G$ being eventually synchronized given the information on $F$. The quantity $  \overline{p}_{T}(G)$ above is the averaged predictive probability that $G$ will eventually synchronize.

	In order to effectively compute the expectation in \eqref{eq:global_predictive_prob}, we use the Monte Carlo approximation along an MCMC trajectory of $k$-paths. That is, our sampling oracle generates a sequence of $k$-paths $(\x_{s})_{s\ge 0}$ in $G$ that forms an irreducible Markov chain. By Lyu et al.\cite[Thm. 2.3]{lyu2023sampling}, we have that almost surely, 
	\begin{align}\label{eq:sample_mean_p}
		\overline{p}_{T}(G) = \lim_{N\rightarrow\infty} \left(p_{T;\, N}(G):= \frac{1}{N}\sum_{s=1}^{N} p_{T}(F[\x_{s}]) \right)
	\end{align}
	where $F[\x_{s}]$ denotes the subgraph of $G$ induced on the nodes in the $s$th $k$-path $\x_{s}$. The above sample average in the right-hand side of \eqref{eq:sample_mean_p} can be computed recursively without storing all past samples. This gives us the following recursive algorithm for computing the approximate predictive probability $p_{T;s}(G)$ for all $s\ge 1$:
	\begin{align}\label{eq:recursive_avg_mcmc}
		\begin{cases}
			\mathcal{X}_{s} &\leftarrow \textup{CAT on subgraph $F[\x_{s}]$} \\
			p^{(s)} &\leftarrow \sigma( \Beta^{T} \MAT(\mathcal{D})^{T}\VEC(\X_{s}) ) \\
			p_{T;\, s}(G) &\leftarrow \left(1-\frac{1}{s}\right)p_{T;\, s-1}(G) + \frac{1}{s} p^{(s)},
		\end{cases}
	\end{align}
	where $\sigma(x)=\frac{\exp(x)}{1+\exp(x)}$  and $(\mathcal{D}, \Beta)$ is a given hyperparameter for LLDM. The recursion \eqref{eq:recursive_avg_mcmc} can be executed over arbitrarily many MCMC samples efficiently. 
	
	Furthermore, by Lyu et al.\cite[Thm. 2.23]{lyu2023sampling}, the recursive averaging \eqref{eq:recursive_avg_mcmc} is guaranteed to converge to the population mean \eqref{eq:global_predictive_prob} exponentially fast in $N$. That is, for each $\delta>0$ and for all $N\ge 1$,
	\begin{align}\label{eq:thm:McDiarmids_1}
		\mathbb{P}\left( \left| \overline{p}_{T}(G) - p_{T;N}(G)  \right| \ge \delta \right) <  2\exp\left( \frac{-2\delta^{2}N}{9 \tau_{mix}} \right),
	\end{align}
	where $\tau_{mix}$ is the mixing time of the standard lazy simple symmetric random walk on $G$, which depends on the size and topology of $G$. In practice, we observe that $p_{T;N}(G)$ converges quickly to $\overline{p}_{T}(G)$ in many problem instances, see Figure \ref{fig:Recursive_averaging_example}\textbf{a}.
	
	\section{Learning hyperparameters}
	
	In this section, we discuss how to learn the model hyperparameters, the regression coefficient vector $\Beta$, and the dictionary of latent dynamics filters $\mathcal{D}$. 
	
	\subsection{Generating the training data set} 
	
	Suppose we have $2m$ observed coupled oscillator dynamics (Kuramoto, FCA, or GHM) $(X^{(j)}_{t})_{0\le t < T}$ for $j=1,\dots,2m$ on the whole graph $G$, and let $y_{j}\in \{0,1\}$ denote the indicator that $X^{(j)}_{t}$ synchronizes as $T\rightarrow\infty$. We assume half the dynamics are synchronizing ($y_{1},\dots,y_{m}=1$) and the other half are non-synchronizing ($y_{i+1},\dots,y_{2m}=0$). Sample $k$-node subgraphs $F_{1},\dots,F_{N_{0}}$ in $G$ (we use $k \in \{10, 15, 20, 25, 30\}$) are sampled through our subgraph sampling oracle described in Section \ref{sec:path_sampling}. We then restrict the first dynamics $(X^{(1)}_{t})_{0\le t < T}$ on the subgraphs $F_{1},\dots,F_{N_{0}}$. This gives $N_{0}$ training examples $(\mathcal{X}_{1},y_{1}),\dots, (\mathcal{X}_{N_{0}}, y_{1})$, where each $\mathcal{X}_{k}$ is the CATs of shape $k\times k\times T$ encoding $(X^{(1)}_{t})_{0\le t < T}$ restricted on $F_{1}$. Next, we restrict the second dynamics $(X^{(2)}_{t})_{0\le t < T}$ on the subgraphs and obtain training examples $(\mathcal{X}_{N_{0}+1},y_{2}),\dots, (\mathcal{X}_{2N_{0}}, y_{2})$, and so on. In total, this creates a training data set consisting of pairs $(\X_{j},y_{j})$, $j=1,\dots,N$ ($N=mN_{0}$), where each $\X_{i}$ is the $k\times k \times T$ CAT of an observed dynamics on subgraph $F_{i}$ and and $y_{i}$ is the corresponding synchronization indicator of the dynamics on the whole graph $G$. 
	
	We also propose an alternative way to generate a training data set when the whole graph $G$ is not available. Suppose we have the same $k$-node subgraphs $F_{1},\dots,F_{N_{0}}$ in $G$ sampled through our subgraph sampling oracle described in Section \ref{sec:path_sampling}, but assume that generating many instances of the coupled oscillator dynamics on the whole graph $G$ is computationally prohibitive. In this case, we can simply run the $2m$ dynamics on the sampled subgraphs $F_{i}$ and record whether the subgraph dynamics synchronize or not with the indicator variable $\tilde{y}_{2m(i-1)+j}$. Denoting by $\widetilde{\X}_{2m(i-1)+j}$ the $k\times k\times T$ CAT of the $j$th dynamics solely run on $F_{i}$, this gives us the training examples $(\widetilde{\X}_{\ell}, \tilde{y}_{\ell})$ for $\ell=1,\dots, N(=m N_{0})$. 
	
	\subsection{How to learn $\beta$ given $\mathcal{D}$}\label{sec: b_given_d}
	
	Once we have $\mathcal{D}$, we can estimate the regression coefficients in $\Beta$ from a set of training examples by solving the standard logistic regression optimization problem. Namely, using the latent dynamics filters in $\mathcal{D}$, we can form the $N\times r$ \textit{proximity score matrix} $\mathbf{H}$, whose $(i,j)$ coefficient is given by 
	\begin{align}\label{eq:def_proximity_mx}
		\mathbf H{[i,j]} := \langle \X_{i}, \F_{j} \rangle.
	\end{align}
	Note that the $i$th column of $\mathbf{H}$ gives the proximity score vector for the $i$th observation $X_{i}$. Then joint log-likelihood of observing $(y_{1},\dots, y_{N})$ under \textup{LLDM} is 
	\begin{align}
		&\log L(y_{1},\dots,y_{n}\,|\, \X_{1},\dots,\X_{n}) 
		= \sum_{i=1}^{n} \log \P(Y=y_i\,|\, \mathbf H)  \\
		\quad &= \sum_{i=1}^{N} \left\{y_i\log \pi_i+(1-y_i)\log(1-\pi_i)\right\},
	\end{align}
	where $\pi_{i}$ is the predictive probability for $\X_{i}$ under the regression parameter $\Beta$ and the known latent dynamics filters in $\mathcal{D}$. We can then estimate the corresponding regression coefficients, $\hat \beta$,
	from the above joint log-likelihood function with maximum likelihood estimation (MLE) as
	\begin{align}
		\hat{\Beta} \in \argmax_{\Beta\in \R^{r}}  \,\, \sum_{i=1}^{N} y_{i} \log \pi_{i} + (1-y_{i}) \log (1-\pi_{i}).
	\end{align}
	The above is an unconstrained convex optimization problem, which can be solved by standard first-order optimization algorithms such as gradient descent \cite{boyd2004convex}. When the two classes (i.e., synchronizing and non-synchronizing) are sufficiently balanced in the training data set, then one can employ faster second-order methods such as Newton-Raphson with numerical stability \cite{fisher1922mathematical}. 
	
	As in generalized linear model theory, we have asymptotic normality of the MLE $\hat{\Beta}$ estimated from independent random samples $(\mathbf{h}_{i}, Y_{i})$ for $i=1,\dots,N$ as $N\rightarrow\infty$. Namely, we assume that for a true model parameter $\Beta^{*}$, 
	\begin{align}
		\textup{logit}(\E[Y_{i}\, |\, \mathbf{h}_{i} ]) = \langle \Beta^{*}, \mathbf{h}_{i} \rangle   + \eps_{i}, 
	\end{align}
	where $\mathbf{h}_{i}$ denotes the $i$th row of $\mathbf{H}$ and $\eps_{i}$s are i.i.d. normal random variables with a constant variance $\sigma^{2}>0$  \cite{nelder1972generalized}. 
	
	The MLE $\hat{\Beta}$ is known to converge to $\Beta^{*}$ with asymptotically normal fluctuation as the sample size tends to infinity. More precisely, 
	\begin{align}
		\frac1\sigma\sqrt{n}(\hat \beta-\beta^{*})\xrightarrow{D}N(\mathbf 0,\Sigma^{-1}),
	\end{align}
	given that the sample covariance matrix $\frac{1}{N}\mathbf{H}^T\mathbf{H}$ converges to a limiting covariance matrix $\Sigma$ as $N\rightarrow \infty$ \cite{billingsley2013convergence}. This provides a statistically powerful mechanism for determining important covariates, or features, that affect the response according to the LLDM, which we demonstrate in Section \ref{sec:goodness-of-fit}. In the following sections, we propose how to compute the dictionary $\mathcal{D}$ of latent dynamics filters. 
	
	A detailed section about the training data can be found in Appendix \ref{app:training_details}.
	
	\subsection{How to learn $\mathcal{D}$ from observed dynamics $\X_{1},\dots,\X_{N}$}\label{sec:learn_filters}
	
	We now know how to learn $\Beta$ given $\mathcal{D}$. In this subsection, we propose some methodologies to learn the dictionary $\mathcal{D}=[\F_{1},\dots,\F_{R}]$ of latent dynamics filters from the observed dynamics in the form of CATs $\X_{1},\dots,\X_{N}$.

	A na\"{i}ve choice of $\mathcal{D}$ is the set of all observed CATs $\X_{i}$ for $i=1,\dots,N$. This means that we regard every single observed dynamic becomes a latent dynamics filter. This choice of $\mathcal{D}$ is undesirable since it is computationally expensive to use a very large number $N$ of filters in our model and also it does not provide any reduced-dimensional representation of the observed dynamics, which hinders the interpretability of our method. Instead, we employ a matrix factorization-based approach (See Appendix \ref{app: nmf}) to learn a small set of bases filters of size $R$, where $R\ll N$ (typically $R\in \{ 2, 8, 25, 100\}$). Principal component analysis (PCA) is a popular tool for extracting key features and reducing the dimensionality of the data set, but it is not suitable for our purpose since we desire non-negative basis elements for the CATs so that we can interpret the basis elements (filters) as representing latent dynamics on latent subgraphs. Hence, we employ \textit{ non-negative matrix factorization} (NMF) (with vectorizing the tensor input) to extract nonnegative latent $k\times k \times T$ latent dynamic filters $\mathcal{F}_{1},\dots,\mathcal{F}_{R}$. See Section \ref{app: nmf} for background and details on NMF.
	\begin{description}
		\item[(i)] \textbf{Feature extraction from observed dynamics:} Use \textit{nonnegative matrix factorization} (NMF) to extract $R$ nonnegative basis tensors $\F_{1},\dots,\F_{R}$ from $\X_{1},\dots,\X_{N}$. 
	\end{description}
	This approach is especially useful when we may not have much prior knowledge about the underlying graph and dynamical system to meaningfully distill out parts of the training data.
	
	There are a number of sufficient conditions on coupled oscillator systems that are guaranteed to lead to global synchronization or non-synchronization (see Sec. \ref{sec:sufficient_conditions}). In scenarios like ours, where one has some existing knowledge about the dynamical system at hand and might want to emphasize higher and finer-detailed interpretability of the learned filters, one can `distill out'  a certain portion of training data before we apply NMF for interpretable feature extraction. This suggests the following variant of the previous method \textbf{(i)} of learning $\mathcal{D}$ from the observed dynamics: 
	\begin{description}
		\item[(ii)] \textbf{Feature extraction from observed dynamics after distillation:} Subsample the `most relevant' observed dynamics $\X_{i_{1}},\dots, \X_{i_{m}}$ using knowledge from coupled oscillator systems, for e.g., dense/sparse underlying graphs and concentrated initial configuration (see Section \ref{sec:sufficient_conditions}) and then apply NMF as in \textbf{(i)}.
	\end{description}
	Note that \textbf{(ii)} is an example of utilizing the tight interplay between knowledge of these coupled oscillators that originates from the literature along with machine learning-based knowledge.

	\subsection{How to learn $\mathcal{D}$ and $\beta$ jointly} 
	
	Combining the procedures in the previous two sections, we can first learn a dictionary of learn latent dynamics filters $\mathcal{D}$ from the observed dynamics $\X_{1},\dots,\X_{N}$ without using the synchronization indicators $y_{1},\dots,y_{N}$ and then learn the regression coefficients in $\Beta$ using the learned $\mathcal{D}$ with synchronization indicators $y_{1},\dots,y_{N}$. However, it is also possible to learn $(\mathcal{D},\Beta)$ jointly from the training examples $(\mathcal{X}_{i},y_{i})$ for $i=1,\dots, N$, and we argue that such joint learning is beneficial for our goal of addressing Problem \ref{problem:main2}. Recall that the latent dynamics filters $\F_{1},\dots,\F_{R}$ in $\mathcal{D}$ should ultimately satisfy the following two objectives: (1) (\textit{ data reconstruction}) They represent patterns in dynamics on subgraphs that are rich enough so that any observed dynamics on a subgraph can be approximately reconstructed by a (nonnegative) linear combination of them; and (2) (\textit{ class-discrimination}) they represent patterns in dynamics on subgraphs that are the most effective in discriminating eventual synchronization and non-synchronization on $G$. If we learn $\mathcal{D}$ only from the observed dynamics as in Section \ref{sec:learn_filters}, it may satisfy (1), but not necessarily (2). 
	
	From this perspective, we also propose to use \textit{supervised matrix factorization} (SMF)\cite{lee2022supervised, lee2023supervised} to learn low-rank latent dynamic factors that offer interpretable, data-reconstructive, and class-discriminative features, addressing the challenge of satisfying the two objectives (1) and (2) simultaneously. SMF is similar to NMF in that it extracts a set number of latent features from the observed data set for interpretable dimension reduction, but the matrix factorization process is supervised by the class labels so that the latent features can also be class-discriminative. For this, we formulate a non-convex constrained optimization problem to jointly learn $\Beta$ and $\mathcal{D}$ and solve that problem iteratively via block-coordinate descent type methods (see Appendices \ref{app:training_details} and \ref{app: SMF}). 
	
	Note that this approach entails a setup for LLDM that is similar to that of a two-layer neural network \cite{rosenblatt1961principles, rumelhart1985learning} with one input, one hidden, and one output layer (See Figure \ref{fig:GLLM_scheme}\textbf{a}). Hence, it is possible to use backpropagation instead of an SMF-based approach to jointly learn the parameters that represent our filters $\F_{1},\dots,\F_{R}$ and regression coefficients $\beta$ by using $\X_{1},\dots,\X_{N}$ as our input layer elements which aim to predict the synchronization indicators $y_{1},\dots,y_{N}$. However, using traditional backpropagation for training a feedforward neural network would result in filters with significantly degraded interpretability, due to the loss of the nonnegativity constraint. Bassi et al.\cite{bassi2022learning} showed that complicated training regimes involving deep feedforward neural networks (FFNN) and long-term recurrent convolutional networks (LRCNs) \cite{donahue2015long} perform quite similarly to other simple algorithms like random forest, gradient boost, and logistic regression on the task of predicting synchronization. We further demonstrate that in Section \ref{sec:results}, where the simple and interpretable architecture of LLDM manages to match or outperform the performance of a multilayer neural network architecture that is much more complex in comparison, while also maintaining the interpretability of the learned filters.
	
	\section{Results} \label{sec:results}
	
	We now report on the performance of \textup{LLDM} for synchronization prediction tasks. We consider two cases, where (1) one seeks to predict the synchronization indicator on small subgraphs using the dynamics on the subgraphs observed during a short time period; and (2) one seeks to predict the synchronization indicator on the whole (parent) graph using the dynamics restricted on several subgraphs observed during a short time period. We also discuss the goodness-of-fit of our model, as well as the interpretability of the learned filters. Experimental details, hyperparameter choices, and model architectures can be found in Appendix \ref{app:training_details}.
	
	\subsection{Networks}
	
	In our experiments, we take the large-scale graph $G$ to be one of the following three types of networks, (1) \dataset{UCLA}, (2) \dataset{Caltech}, and (3) networks generated from the Newmann-Watts-Strogatz (\dataset{NWS}) model. \dataset{UCLA} and \dataset{Caltech} networks are part of the \texttt{FACEBOOK100} dataset \cite{traud2012social}, where the nodes represent users in the respective Facebook networks, and the edges encode Facebook `friendships' between these accounts. Furthermore, for our third large network, we generate a single connected graph using the Newman–Watts–Strogatz (\dataset{NWS})  \cite{newman1999renormalization} small-world network model with $n = 20000$ nodes, $k = 1000$ nearest neighbors in the circulant initial graph, and each non-adjacent pair of nodes gets a new edge independently with probability $p = 0.5$. Lastly, we also generate 500 instances of \dataset{NWS} networks with $n=300$ nodes, $k=12$ nearest neighbors, and shortcut edge probability $p=0.4$. The basic summary statistics of all these three networks can be found in Table \ref{tab:networks}.
	
	\begin{table}[!ht]
		\caption{\label{tab:networks}
			Basic Graph Statistics of the Large-scale Graph Networks.
		}
		\begin{ruledtabular}
			\begin{tabular}{c|cccc}
				\textbf{Networks} & \dataset{UCLA} & \dataset{Caltech} & \dataset{NWS} & \dataset{NWS}'\\
				\hline
				$\#$ of graphs & 1 & 1 & 1 & {500} \\
				$\#$ of nodes & 20467 & 769 & 20000 & 300\\
				$\#$ of edges & 747613 & 16656 & 1.49e+7
				& {2.52e+3 $\pm$ 19.4}\\
				Edge density & 0.0036 & 0.0564 & 0.0750
				& {0.0562}  \\
			\end{tabular}
		\end{ruledtabular}
	\end{table}

	\subsection{Sufficient conditions for synchronization}\label{sec:sufficient_conditions}
	
	Some properties of the underlying graphs or dynamics themselves are well known to be critical for the synchronization behavior of a system of coupled oscillators. We state some of these well-known properties below, which come from the traditional literature on coupled oscillator systems. These conditions will be the basis of our baseline predictor, whose performance gives a sense of how easy or difficult a given synchronization prediction problem is. 
	
	It is well-known that coupled oscillator systems on dense graphs are relatively easier to synchronize compared to systems on sparse graphs. For example, Kassabov, Strogatz, and Townsend recently showed that Kuramoto oscillators with identical natural frequency on a connected graph where each node is connected to at least 3/4$^{\textup{th}}$ of all nodes are globally synchronizing for almost all initial configurations \cite{kassabov2021sufficiently}. Moreover, a more intricate analysis of Kuramoto oscillators on networks shows that it is possible to generate dense circulant networks that are just able to prevent global synchronization, and sparse circulant networks which tend to globally synchronize \cite{townsend2020dense}.
	
	For FCA, it is known that the dynamics synchronize on a path for a $\kappa-$coloring configuration if $\kappa \ge 3$ and further that for finite trees with maximum degree $\Delta$, FCA dynamics do not synchronize if $\Delta \ge \kappa$ and synchronize if $\Delta < \kappa \le 6$\cite{lyu2015synchronization, lyu2016synchronization}. Also, GHM tends to not synchronize on complete or highly dense graphs. \cite{bassi2022learning} also showed that GHM synchronizes on paths. 
	
	The next sufficient condition is the \textit{concentration principle}, which is a fundamental observation in coupled oscillators and has been widely used in the clock synchronization literature \cite{nishimura2011robust, klinglmayr2012guaranteeing, lyu2018global} as well as multi-agent consensus problems \cite{moreau2005stability, chazelle2011total}. This principle, stated below, follows from the fact that the phase difference between any two nodes, when isolated from the rest, monotonically decreases to zero, assuming an identical oscillation frequency.  
	
	\begin{description}
		\item[Concentration Principle] \textit{Consider an arbitrarily connected graph $G.$ For Firefly Cellular Automata (FCA), as well as the Kuramoto Model (KM) with identical intrinsic frequency, the given dynamics on $G$ synchronize if all phases at any given time are confined in an open half-circle in the phase space $\Omega$. Furthermore, if all states used in the configuration $X_t$ are confined in an open half-circle for any $t \geq 1$, then the trajectory on $G$ will eventually synchronize.}
	\end{description}
	
	An open half-circle refers to any arc of length $< ~\pi$ for the continuous phase space $\Omega = \mathbb{R}/2\pi\mathbb{Z}$ and any interval of $< \kappa/2$ consecutive integers (mod $\kappa$) for the discrete phase space $\Omega = \mathbb{Z}/\kappa\mathbb{Z}$. This confinement in an open half-circle is what we define a phase configuration as being `concentrated'. Further, since the concentration principle does not hold for the Greenberg-Hastings model, we define a phase configuration for GHM to be `concentrated' if it is synchronized. The baseline predictor we use for our experimental validation in Section \ref{sec:results} is based on the concentration principle: \textit{Predict synchronization is the phase configuration at any time during the observed dynamics is concentrated; otherwise flip a fair coin.} See Appendix \ref{app:training_details} for details. 
	
	We take advantage of some of these sufficient conditions when we formulate a theory-informed data distillation approach in Section \ref{sec:learn_filters}
	to learn latent dynamic filters for \textup{LLDM}. 
	
	\subsection{Model validation I: Subgraph level}\label{sec:res_local}
	
	Here we discuss the performance of the latent linear dynamics model (\textup{LLDM}) at the subgraph level for synchronization prediction. We discuss the results of synchronization prediction using \textup{LLDM} based on joint-optimization using $\textup{SMF}$ (labeled ``LLDM'' in Table \ref{tab:acc_table}) as well as the theory-informed data distillation approach (see Sec. \ref{sec:learn_filters} \textbf{(ii)}) (labeled ``\textup{LLDM-T}'' in Table \ref{tab:acc_table}). We consider the case of $10$-, $20$-, and $30$-node subgraphs sampled using the subgraph sampling oracle in Sec. \ref{sec:path_sampling} from the 20000-node \dataset{NWS} parent graph (in Table \ref{tab:networks}). We generated training and testing data sets by running Kuramoto, FCA, and the GHM dynamics on the sampled subgraphs.
	See Appendix \ref{app:training_details} for more details on generating data sets for the experiments.

	\begin{table}[!ht]
		\caption{\label{tab:acc_table}
			Prediction accuracy of various methods for FCA and Kuramoto dynamics on subgraphs with $k$ number of nodes where $k \in \{10, 20, 30\}$ sampled from a large-connected \dataset{NWS} parent graph. All accuracy values are an average of $5$ seeds. The highest accuracy for each setting is indicated in \textbf{bold} font.} 
		\begin{ruledtabular}
			\begin{tabular}{c|ccc|cccc}
				&&FCA&
				&&Kuramoto&\\
				\hline
				&$k=10$&$k=20$&$k=30$
				&$k=10$&$k=20$&$k=30$\\
				\hline
				Baseline&80.2&69.6&69.4&56.5&66.1&64.3 \\
				LogReg&92.7&94.8&76.4&83.5&80.1&80.2 \\
				FFNN&92.2&\textbf{95.4}&82.8&\textbf{84.2}&\textbf{83.4}&79.2 \\
				\hline
				LLDM-T (R=2)&82.7&81.3&75.9&74.3&76.0&80.1\\
				LLDM-T (R=8)&86.5&81.9&75.7&75.2&77.1&\textbf{80.3}\\
				LLDM (R=2)&91.5&93.2&77.2&81.4&77.5&78.3 \\
				LLDM (R=8)&\textbf{93.0}&94.2&\textbf{84.8}&82.1&78.4&78.4 \\
			\end{tabular}
		\end{ruledtabular}
	\end{table}
	
	We used extremely low-rank LLDM with $R\in \{2,8\}$ and compared the prediction accuracy with that of the baseline predictor, logistic regression, and feedforward neural network. A subset of our experimental results is reported in Table \ref{tab:acc_table}. For a detailed set of results on different subgraph sizes, parent networks, and dynamics models, see Appendix \ref{app:extended_subgraph_results}. 
	
	Consider the performance of LLDM using SMF on FCA and Kuramoto dynamics in Table \ref{tab:acc_table}. We observe that LLDM outperforms the baseline in all settings. We also observe that LLDM performs well across all three subgraph sizes for both the Kuramoto and FCA dynamics. In some cases, our method outperforms logistic regression as well as FFNN, while in other cases, LLDM is still quite close to the prediction accuracies of black-box methods.

	Additionally, there seems to be a trade-off between choosing the number $R$ of dictionary atoms to learn from our data. Generally, a higher rank parameter $R$ leads to a higher accuracy as more atoms can effectively capture more fine-grained features of the data. Furthermore, it can be seen that for the problem instances we created, it is easier to predict the synchronization of the FCA dynamics than for Kuramoto dynamics. Furthermore, the prediction task becomes harder as the subgraph size increases, as perhaps the synchronization indicator on a larger subgraph depends on more complex patterns in dynamics. We also observe that the theory-informed LLDM-T in most settings performs worse compared to the SMF-based LLDM, but the accuracies are in general still comparable. This shows the potential of harnessing a theory-informed approach if we have appropriate knowledge about the dynamical systems, especially in settings where we have a large number of noisy observations. 
	
	\subsection{Model validation II: Global level}\label{sec:res_global}
	
	Next, we use LLDM to predict synchronization on a large graph using information only at the subgraph level. Here, we take $G$ to be one of the 500 instances of the 300-node graphs in the dataset \dataset{NWS}'  in Table \ref{tab:networks}. We run the FCA and Kuramoto dynamics on $G$ for $T'=50$ and 100 iterations, respectively, and let $y_{G}$ be the indicator of whether the system on $G$ is globally synchronized at time $T'$. Half of them globally synchronize at time $T'$ (so $y_{G}=1$) and the other half do not (so $y_{G}=0$), which are split into 80$\%$ training set and 20$\%$ testing set.  From each $G$, we sample a single trajectory of 50 iterations of the MCMC $k$-path ($k \in \{10, 20, 30\})$ motif sampling algorithm (see Section \ref{sec:path_sampling}), which gives a sequence of $k$-node subgraphs $F_{1},\dots,F_{50}$. We then restrict the dynamics on $G$ on the sampled subgraphs, thereby creating pairs $(\X_{1},y_{G}),\dots,(\X_{50}, y_{G})$ of $k\times k\times T$ CATs and (global) synchronization indicator, where $T$ is varied between 10 and 100. In this way, we create a total of 500*50 pairs of observed subgraph dynamics and synchronization indicators.
	
\begin{figure}[!ht]
		\centering
		\includegraphics[width=0.48\textwidth]{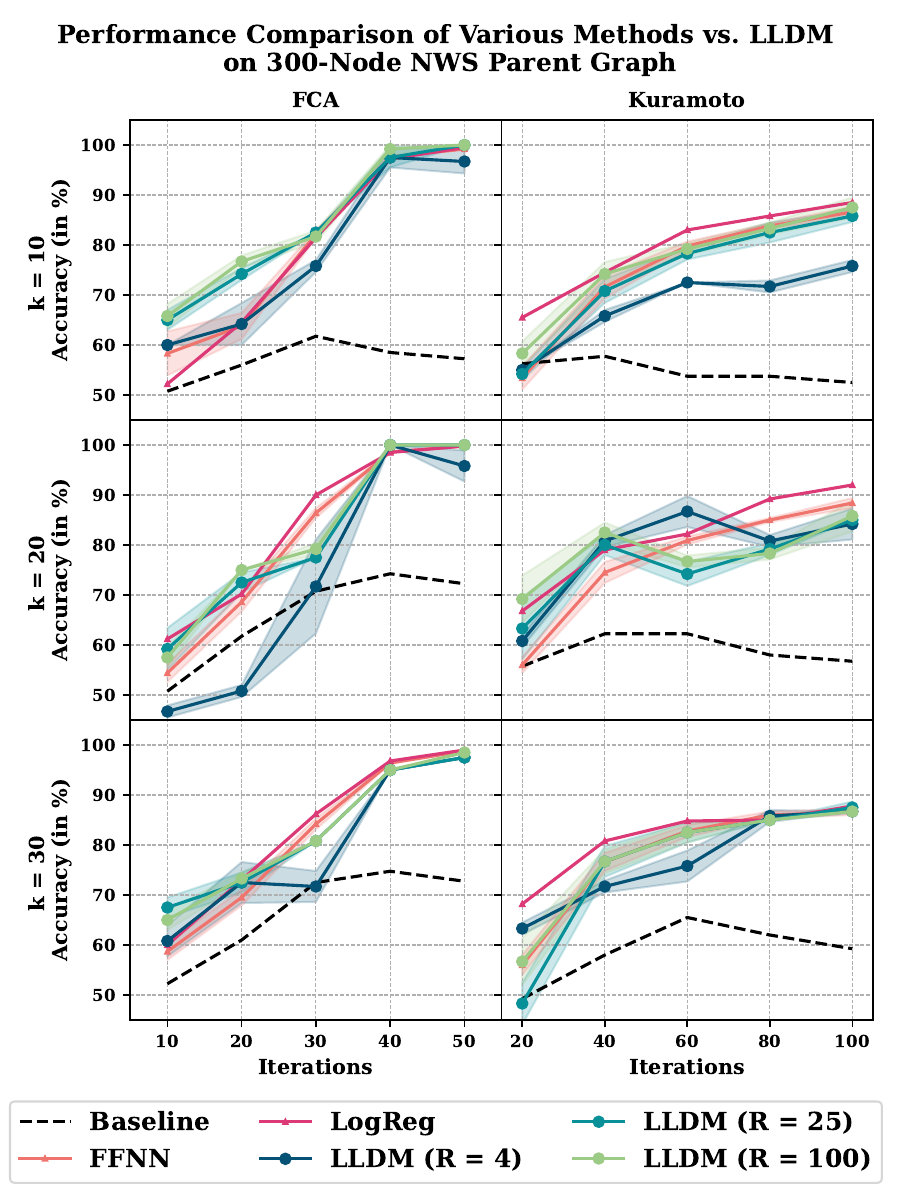}
		\caption{Prediction accuracy of various methods on large 300-node \dataset{NWS} graphs by using subgraphs of sizes 10, 20, and 30. The prediction accuracies and the standard deviation error shades are shown for varying levels of iterations of dynamics data shown in these models.} 
		\label{fig:LLDM_global}
	\end{figure}

	We use a block-minimization algorithm for SMF \cite{lee2023supervised} on the training set to jointly learn the dictionary of latent dynamic filters $\mathcal{D}$ and the vector of regression coefficients $\Beta$, where we used three rank parameters $R\in \{4,25,100\}$. We can then use the trained LLDM on the testing data set and approximately compute the predictive probability $\overline{p}_{T}(G)$ that the dynamics on $G$ will eventually synchronize with the sample average in \eqref{eq:sample_mean_p} and \eqref{eq:recursive_avg_mcmc}. For more details of the experimental setup, see Appendix \ref{app:training_details}. The experimental results for predicting the synchronization on $G$ from the subgraph dynamics are shown in Figure \ref{fig:LLDM_global}. 
	
	It is important to note that our method of aggregating subgraph-level predictions to form a single global-level prediction, as we proposed in Section \ref{sec:subgraph_to_parent}, is a general procedure that can be applied to any subgraph-level predictor.  That is, if we have a model that computes a predictive probability $p_{T}(F)$ of global synchronization on $G$ based on a $T$-iterations of dynamics observed on a subgraph $F$, then we can compute the predictive probability for many subgraphs $F$ and take their mean to be the predictive probability of global synchronization as in \eqref{eq:global_predictive_prob}. We compare the performance of LLDM for rank $R\in \{4,25,100\}$ against the baseline predictor, logistic regression, and FFNN, where we use the same local-global aggregation \eqref{eq:global_predictive_prob} for all methods for a fair comparison. 
	
	We first observe that all three methods, logistic regression, FFNN, and LLDM outperform the baseline for both dynamics models on almost all subgraph sizes. Furthermore, it holds consistent from the subgraph-level experiments that prediction of synchronization on FCA is in general easier than that of Kuramoto, as we see that most methods can perfectly predict the emergent properties of the parent graph when shown all $50$ iterations of the evolution of dynamics, while for Kuramoto the accuracy peaks at around $90\%$ when the model is shown $100$ iterations of dynamics evolution.
	
	Second, the global-level prediction task here becomes easier if we increase the subgraph size, contrary to the subgraph-level prediction task in Section \ref{sec:res_local}, as more information on the global dynamics is revealed by observing dynamics on larger subgraphs.

	\begin{figure}[!ht]
		\centering
		\includegraphics[width=0.45\textwidth]{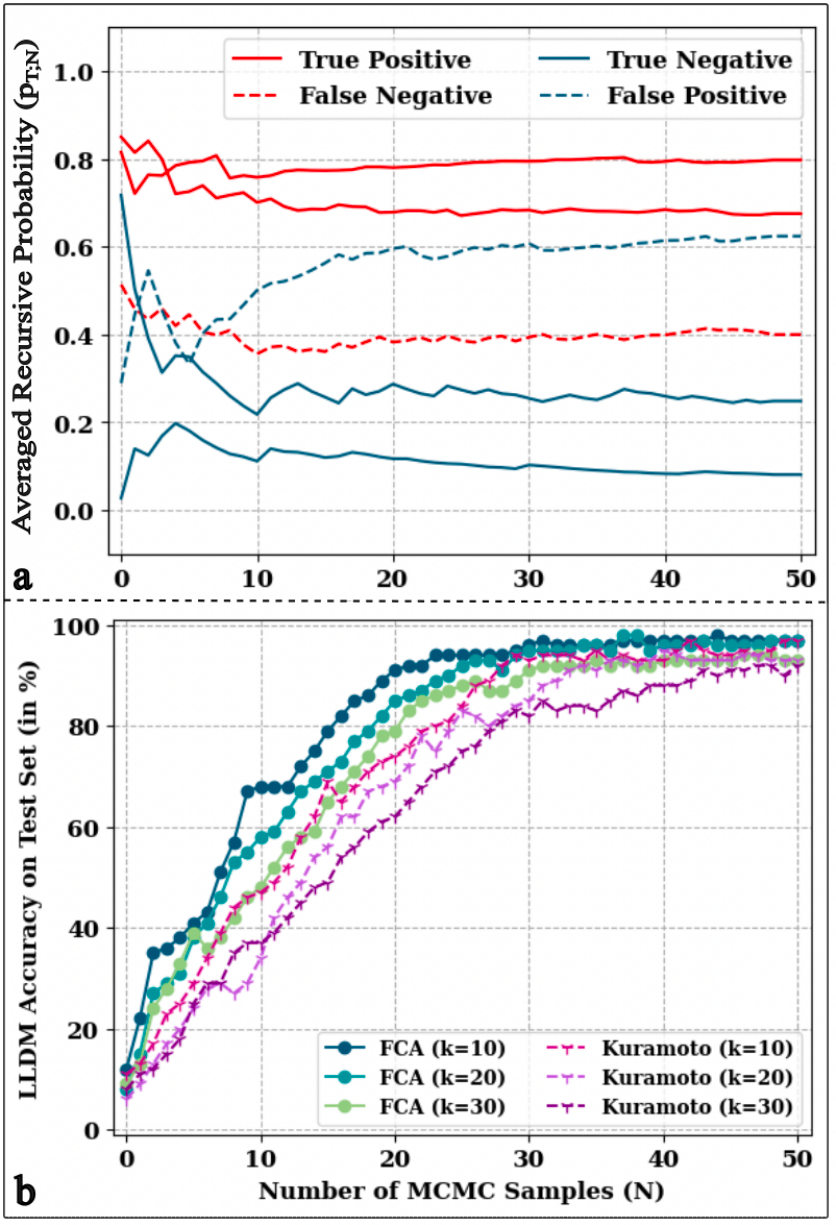}
		\caption{(Panel \textbf{a}) Examples of recursively averaged probability for Kuramoto dynamics on $20$-node subgraphs. We observe that the predicted probability tends to converge in the first 20-25 MCMC samples. (Panel \textbf{b}) Classification accuracy of a pre-trained rank-$4$ LLDM dictionary on the test set. 
		} 
		\label{fig:Recursive_averaging_example}
	\end{figure}

	Third, observe that LLDM with ranks $R=25$ and $R=100$ consistently matches or outperforms all other methods. LLDM with $R=4$ seems to lag behind LLDM with higher ranks and other methods, but still provides a competitive prediction accuracy, indicating that even with just $4$ latent dynamics filters extracted from the vast amount of data, we can predict reasonably and in some cases even match LLDM with much higher ranks and other methods.
	
	Fourth, the error bars in Figure \ref{fig:LLDM_global} represent one standard deviation of the prediction accuracy computed over ten runs of learning $(\mathcal{D},\Beta)$ by SMF and the prediction of LLDM, where the randomness comes from the random initialization of SMF algorithm. From these error bars, we observe that LLDM with rank $R=4$ has a higher fluctuation of prediction accuracy compared to its higher-rank counterparts. This suggests a potential trade-off that comes with selecting the appropriate rank for such large-scale synchronization prediction. While low-rank LLDM proves to be computationally efficient due to being light-weight by learning fewer atoms in the dictionary, this comes at the cost of a higher uncertainty in the prediction due to not a large amount of variation of the data being explained by just $4$ atoms.
	
	Lastly, we investigate how the overall prediction accuracy is affected by the number $N$ of subgraphs we average the predictive probabilities over (by using \eqref{eq:recursive_avg_mcmc}) to get an estimate of the population average \eqref{eq:global_predictive_prob}. Figure \ref{fig:Recursive_averaging_example}\textbf{a} verifies that the recursively averaged predictive probabilities $p_{T;N}$ in \eqref{eq:recursive_avg_mcmc} indeed converge as we increase the sample size $N$. Figure \ref{fig:Recursive_averaging_example}\textbf{b} shows that the overall prediction accuracy on the test set increases linearly as $N$ increases, saturating around $N=25$ samples for FCA and $N=35$ samples for Kuramoto dynamics. 
	
	Overall, these results strongly suggest that our approach of utilizing subgraph-level predictions to extrapolate to the parent graph level helps to tackle the prediction of complex dynamics on large graphs very effectively and in a compute-efficient manner with a simple framework.

	\subsection{Model Validation III: Goodness of Fit -- Linearizing Nonlinear Dynamics}\label{sec:goodness-of-fit}
	
	Recall that LLDM seeks to model the synchronization indicator of nonlinear dynamical systems via a linear representation of some latent dynamical patterns. The successful experiments in the previous sections demonstrate that the synchronization indicator may indeed be modeled as a linear function of some latent features observed in nonlinear dynamics. In this section, we provide further evidence for this claim by utilizing statistical analysis for generalized linear models such as goodness-of-fit and deviance residual plots. 
	
	We first discuss a visual heuristic for goodness-of-fit in linear regression. Let $Y$ denote a univariate response variable and $X=[X_{1},\dots,X_{R}]$ denote a vector of covariates. We assume the conditional expectation $y:=\E[Y\,|\, X]$ of $Y$ given $X$ as a linear function $\hat{y}:=\beta^{T}X$, where $\Beta\in \R^{R}$ is a vector of regression coefficients. However, there could be a higher order dependence between $y$ and the covariates. To see whether a linear model has a goodness-of-fit, we investigate the residual $\hat{r}$ below as a function of the covariates:
	\begin{align}
		\hat{r}:=y- \hat{y}=O(X^p),
	\end{align}
	where $O(X^p)$ denotes a higher-order, polynomial relationship among the covariates in the data. We can plot the residual $\hat r$ by the estimate of the fitted values in regression, $\hat{y}$,
	and in practice, it is often the case that we will see a nonlinear pattern before any model tuning or data transformation is performed in linear regression.
	
	The importance of observing the relationship between residuals by their fitted values is that it serves as a heuristic measure for the goodness-of-fit of a specified linear model. If any nonlinearities were present substantially from this visual heuristic, it would suggest to a practitioner that higher-order interactions between covariates must be incorporated. For example, a cubic residual plot would suggest incorporating cubic order terms between covariates in a more complex model to capture the cubic nonlinearities that may be missed in a simpler linear model of purely first-order covariates. 
	
	In our setting for logistic regression, we utilize `deviance residuals'\cite{mccllagh1989generalized}, which is the appropriate choice for modeling residuals in logistic regression models that are as follows:
	\[d_i:=\text{sign}(e_i)\left[-2(y_i\log \hat p_i+(1-y_i)\log (1-\hat p_i))\right]^{1/2},\]
	for data point $i$, with 
	$e_i=y_i-\hat p_i,$ where $y_i$ is the corresponding label for observation $i$, and $\hat p_i$ is the predicted response after fitting a logistic regression model as in 
	\begin{align}
		\hat p_i:=\frac{\exp( \beta^{T} \mathbf X_{i,:})}{1+\exp( \beta^{T}\mathbf X_{i,:})},
	\end{align}
	where $\mathbf X_{i,:}$ is the $i$th row of $\mathbf X$, corresponding to the $i$th observation.

	\begin{figure}[!ht]
		\centering
		\includegraphics[width=0.48\textwidth]{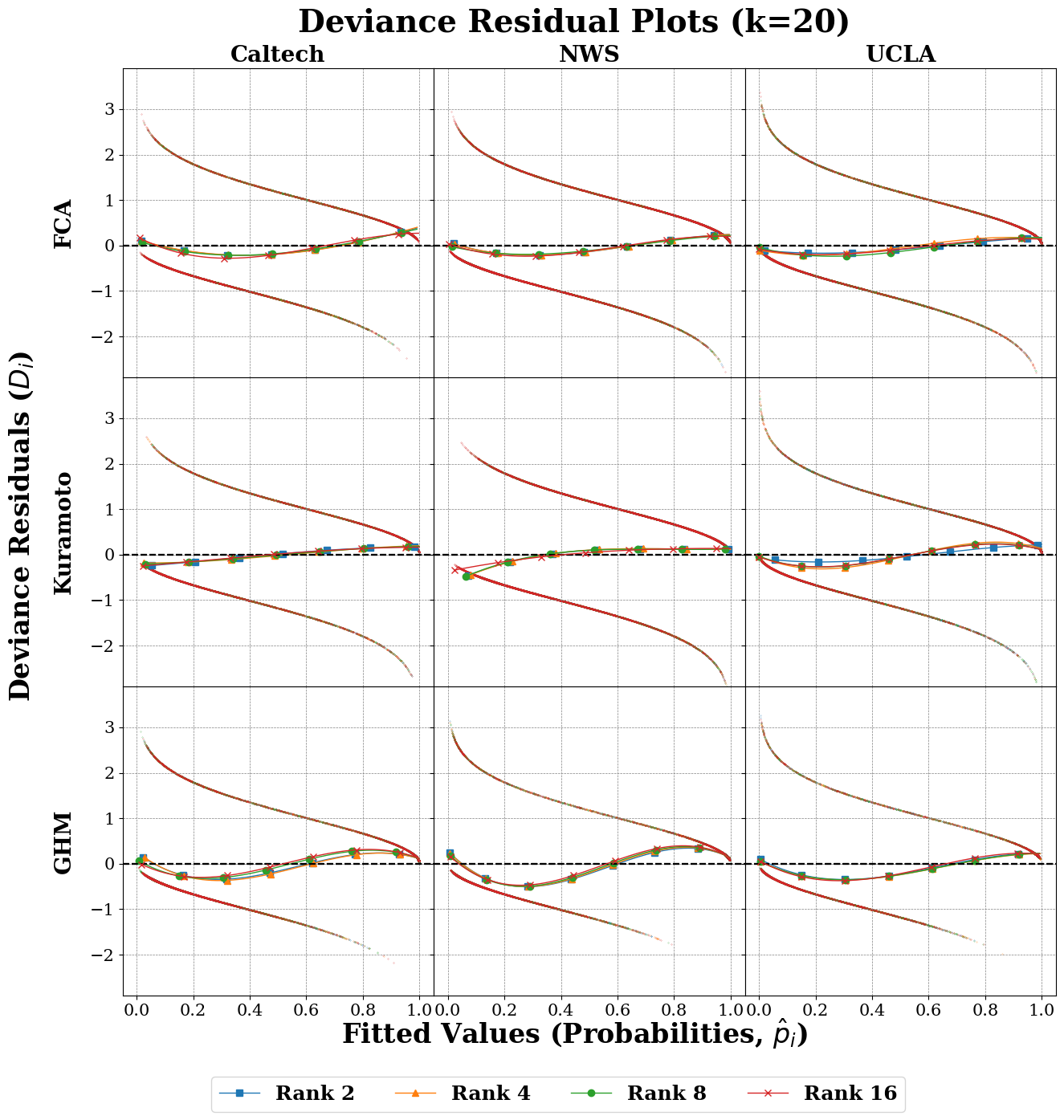}
		\caption{Goodness-of-fit with deviance residuals for LLDM on subgraphs with $k=30$ nodes sampled from \dataset{Caltech}, \dataset{UCLA}, and \dataset{NWS} graphs with FCA, Kuramoto, and GHM models on rank $R\in \{2,4,8,16\}$. We generally observe that across all ranks $R$, our model for 20-node subgraphs has a very good linear fit.
			The univariate smoothing spline for each rank smooths and interpolates the $(\hat p_i, d_i)$ points into a polynomial curve, where we can choose the degree of the polynomial, $k$, as well as a smoothing factor, $s$. We choose $k=3$ and $s=R-\sqrt{2R}$ as a common heuristic choice.
		}
		\label{fig:dev30}
	\end{figure} 
	
	We form a deviance residual plot with pairs of data points between fitted values and deviance residuals, $(\hat p_i, d_i)$, and
	of a smoothing spline approximating a residual curve of $(\hat p_i, d_i)$. 
	(This residual curve can be produced from various options of smoothing spline algorithms.) Since we are interested in modeling the residual data after fitting our LLDM, we utilize a simple univariate spline algorithm
	\cite{dierckx1975algorithm}. 
	The resulting smoothing spline represents the distribution of the deviance residuals away from $y=0$ across possible fitted values such that a more sporadic nonzero curve---especially at the tail ends---would indicate poor model fit, and a flatter, more zero-valued curve would mean sufficient linear fit as the deviance residual appears to be approximately zero across all possible fitted values. 
	
	Figure \ref{fig:dev30} shows multiple goodness-of-fit plots of deviance residuals by fitted values where we have smoothing mostly flat spline curves. In particular, there is a good linear fit for LLDM utilizing ranks 2, 4, 8, and 16 for FCA, Kuramoto, and GHM on all three networks,
	\dataset{UCLA}, \dataset{Caltech} and \dataset{NWS} for subgraphs
	of 30 nodes. More plots are given in the Appendix \ref{app:gof_extended} in Figure \ref{fig:dev15} for subgraphs of various sizes $k=10$, 15, 25, and 30 (we see that there is goodness-of-fit for these models as well). Therefore, LLDM not only performs well in predicting the synchronization of coupled oscillator systems in terms of accuracy but also in representing the synchronization indicator in a latent \textit{linear} form. 
	
	An important observation in Figure \ref{fig:dev30} is that for higher filter matrix ranks there is less linear LLDM fit. An intuitive explanation for this behavior is that the LLDM does not benefit from greater linear model complexity and that additional covariates contribute noise to the model. We see cases where a rank 2 LLDM is more than sufficient to linearly capture the synchronization indicator as a function of the proximity scores of latent dynamic patterns. 
	
	\begin{figure*}[!ht]
		\centering
		\includegraphics[width=1\linewidth]{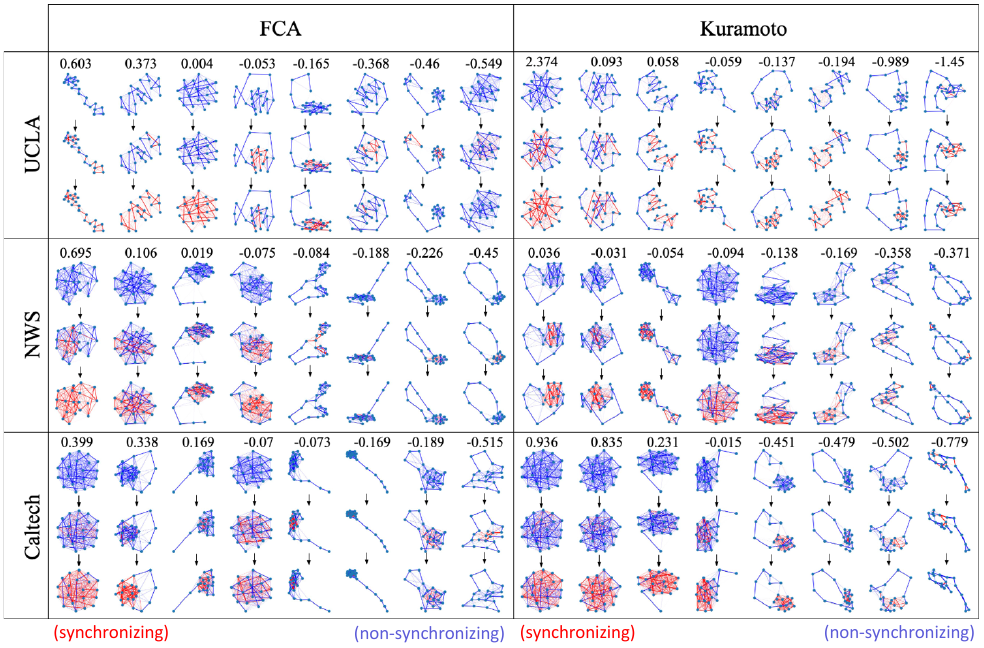}
		\caption{The latent dynamics filters learned by rank-8 SMF for FCA and Kuramoto dynamics on 20-node \dataset{Caltech}, \dataset{UCLA}, and \dataset{NWS} networks. Logistic regression coefficients associated with the latent dynamic filters are shown on top of each filter ($+$ = synchronization and $-$ = non-synchronization). Each latent dynamic filter is a tensor of shape $k\times k\times T$ with $k=20$, and we only show three temporal snapshots of such tensor in each column. The colors on the edges range from blue (indicating a large phase difference) to red (indicating a small phase difference). For instance, the leftmost latent dynamic filter for \dataset{Caltech} subgraphs and FCA dynamics (bottom left) starts with densely connected blue edges ending up with all red edges, indicating a dynamic pattern on densely connected subgraphs with large mutual phase differences leading to synchrony. 
		}
		\label{fig:SMF_filters_fcakura}
	\end{figure*}
	
	\begin{figure}[!h]
		\centering
		\includegraphics[width=1\linewidth]{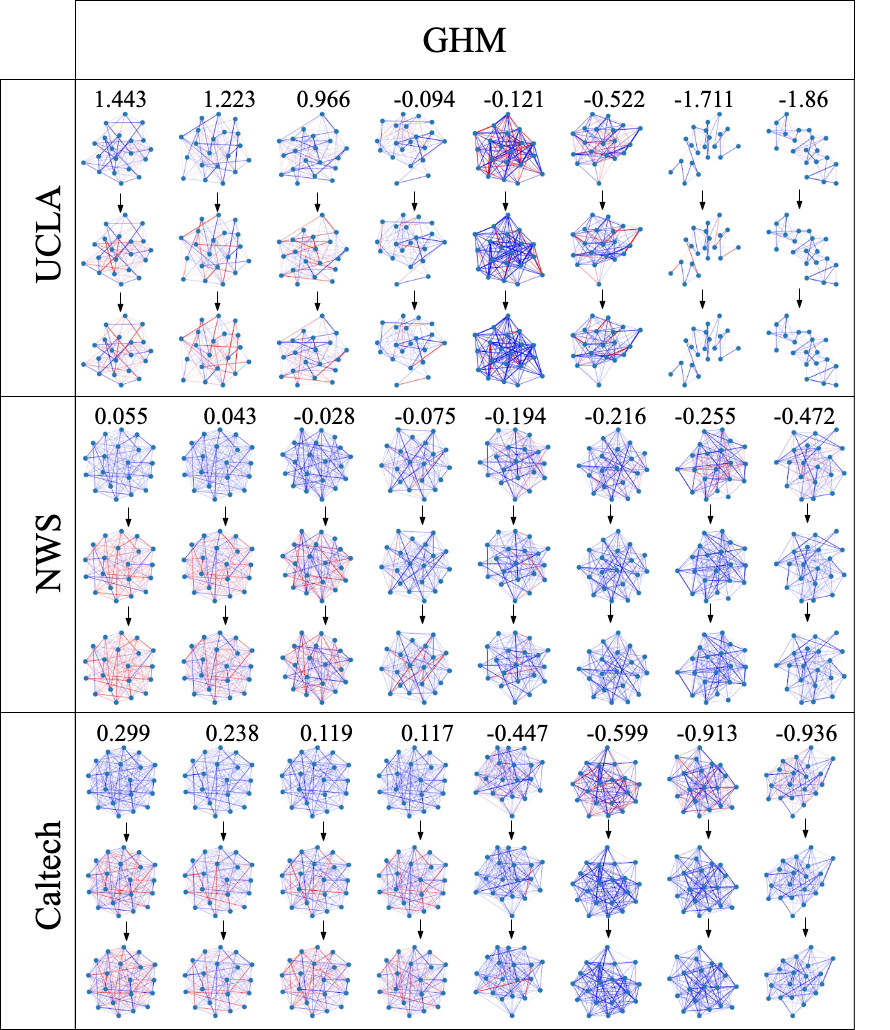}
		\caption{The latent dynamics filters learned by rank-8 SMF for GHM dynamics on 20-node \dataset{Caltech}, \dataset{UCLA}, and \dataset{NWS} networks.}
		\label{fig:SMF_filters_ghm}
	\end{figure}

	\subsection{Interpretability of the learned latent dynamics filters}\label{sec:res_interpret}

	Now that we have validated the performance of LLDM at both the subgraph and global level using accuracy as well as the goodness-of-fit for the coefficients, we now proceed to discuss the \textit{interpretability} of latent dynamics filters (`filters' hereafter), as seen in Figure \ref{fig:SMF_filters_fcakura}, Figure \ref{fig:SMF_filters_ghm}, and in Appendix \ref{appendix:dictionary-plots}. Here, we consider rank $R=8$ filters for LLDM learned by SMF algorithm from the Kuramoto, FCA, and GHM dynamics restricted on the 20-node subgraphs sampled from parent networks of \dataset{Caltech}, \dataset{UCLA}, and \dataset{NWS} in Table \ref{tab:networks}. Similar plots for subgraphs of size $k\in \{10,15,25,30\}$ can be found in Appendix \ref{appendix:dictionary-plots}.

	For each of the tile in both Figures \ref{fig:SMF_filters_fcakura} and \ref{fig:SMF_filters_ghm}, the eight filters are represented by eight columns. Recall that each filter is a CAT of size $k\times k \times T$, representing latent dynamics on $k$-node subgraphs. We only show three equally-spaced snapshots of such tensors as three $k\times k$ matrices in each column, where the time evolution is indicated by the arrows. The colors on the edges range from blue (indicating a large phase difference) to red (indicating a small phase difference). For instance, on the one hand, the leftmost filter for \dataset{Caltech} subgraphs and FCA dynamics (bottom left) starts with densely connected blue edges ending up with all red edges, indicating a dynamic pattern on densely connected subgraphs with large mutual phase differences leading to synchrony. On the other hand, the rightmost filter for \dataset{Caltech} and FCA dynamics ends up with relatively sparse blue edges, indicating that the corresponding subgraph is sparsely connected and the phase differences along its edges are large at time $T$. The filters are in decreasing order for their corresponding logistic regression coefficients (shown on top), where positive (resp. negative) coefficients indicate positive association with eventual (resp. non-)synchronization. So the filters to the left are more representative of what kind of input dynamics and graphs are more likely to synchronize and as we move to the right of the plots, change to being more likely to not synchronize.
	
	We now examine some specific details that can be observed from Figure \ref{fig:SMF_filters_fcakura}, where we consider the case of the FCA and Kuramoto dynamics on \dataset{Caltech} as the underlying graph on the last row of tiles in the figure. We see that for both dynamics, the leftmost filters indicating patterns for synchronization have the densest edges, and almost the entire graph tends to synchronize as time goes on as almost all edges become red. Moving along the columns to the right, we see that almost all filters that have positive or very small  negative coefficients tend to have graphs that are either densely connected or are mostly dense with short paths that go out of a community (a densely connected subset of nodes) structure. Moreover, for such cases, we again observe that the communities that we observe tend to synchronize by the end of the training iteration $T$ but the nodes on the short extending paths do not. We further move right to observe filters that are more likely to not synchronize, observing most of the weakly non-synchronizing filters have a small community structure that weakly synchronizes (see partially communities consisting of red edges) at later times. However, there is a large extending path or cycle that is non-synchronizing. 
	
	Lastly, observe that the filters associated with the most negative coefficients show that the graph structures are very sparse and do not contain any substantial communities. Moreover, almost the entire graph remains not synchronized at the end of the training iteration $T$. Therefore, the latent dynamic patterns captured by these filters with the associated regression coefficients are consistent with the existing knowledge of sufficient conditions for synchronization and non-synchronization. A similar observation holds for the filters learned from the dynamics on subgraphs of \dataset{NWS}. 
	
	Next, we discuss filters learned from oscillator dynamics on subgraphs of \dataset{UCLA}. Recall that \dataset{UCLA} is an order of magnitude sparser than other networks \dataset{Caltech} and \dataset{NWS} (see Table \ref{tab:networks}), so the filters learned from the subgraphs of \dataset{UCLA} show sparser edge density than the ones learned from subgraphs of the other networks. As before, the regression coefficients for Kuramoto dynamics tend to be positively correlated with edge density in the filters, indicating that Kuramoto dynamics on dense subgraphs tend to be synchronizing. The filter with the largest regression coefficient is dense on the whole, and the ones that show weak synchronization are weakly dense yet seem to have a loose community structure and no long-extending paths. Similarly, for the non-synchronizing filters, we observe a very loose to no community structure along with long paths and cycles. However, we see different trends for the filters for the FCA dynamics. There, the regression coefficient seems to depend on the specific topology of the subgraphs and the dynamics on them, rather than just the edge density of the subgraphs we observe the dynamics on. 
	
	Finally, we move on to discuss Figure \ref{fig:SMF_filters_ghm} where we look at the rank $R=8$ filters learned for the GHM dynamics on subgraphs of the three networks. The filters, in this case, look quite similar to each other, unlike the Kuramoto and FCA dynamics in Figure \ref{fig:SMF_filters_fcakura}, which suggests that for GHM dynamics, fewer than $R=8$ filters would suffice to fit our LLDM. Indeed, we obtain similar prediction accuracy for GHM dynamics with $R=2$ as with $R=8$, see Table \ref{tab:acc_table_nws} in Appendix \ref{app:extended_subgraph_results}. 
	
	Furthermore, we observe that most of the filters with a positive association with synchronization seem to be dense overall, and the edges in the filters tend to become red (synchronized) overall within the first half of the training iterations $T$. On the other hand, the filters that correspond to non-synchronization (rightmost columns in Fig. \ref{fig:SMF_filters_ghm}) seem to be relatively sparser or have a path-like structure. There we observe synchronized edges (in red) developing within the filters in time, but most edges remain non-synchronized (in blue).

	\section{Conclusion}
	In this paper, we propose a latent linear model that effectively linearizes highly nonlinear dynamics resulting from coupled oscillators interacting on graphs. A fundamental concept we introduce is `latent dynamic filters', which encode some key dynamical patterns associated with synchronization/non-synchronization of the system, which enable subgraph-level synchronization prediction by using a latent linear model. These filters are directly learned from the data by incorporating supervised matrix factorization techniques. The predictive probability of global synchronization is computed by averaging many subgraph-level predictions along an MCMC subgraph sampling trajectory. 
	
	Our framework has the benefit of being simple and lightweight, and we carried out an extensive study to show that it matches or outperforms traditional black-box methods on the prediction of synchronization of the coupled oscillator dynamics on large graphs. We provide an efficient recursive averaging algorithm to combine many subgraph-level predictions, whose convergence is both theoretically and experimentally justified. A statistical validation of our model was also provided. Finally, we provide a computational framework to extract key patterns responsible for synchronization/non-synchronization from many instances of observed coupled oscillator dynamics. This provides added interpretability of our method for a better understanding of rich nonlinear dynamics of coupled oscillators. We hope that our work will inspire a new line of research harnessing the potential of our simple and interpretable computational framework to help better understand various complex dynamics systems beyond coupled oscillators.
	
	\section{Materials}
	The code for the algorithms and simulations used in this work is provided in \url{https://github.com/zwu363/Interpretable-ML-Sync}. \\
	
	\section*{Acknowledgements}
	
	This work is supported by NSF Grants DMS-2010035 and DMS-2206296.
	
	\section*{References}
	
	\bibliography{mybib}
	
	\clearpage
	\appendix

	\section{Models of Coupled Oscillators} \label{app: oscillators}
	
	Systems of coupled oscillators have been studied over multiple decades \cite{strogatz2000kuramoto,acebron2005kuramoto,nair2007stable,pagliari2011scalable,dorfler2011critical,dorfler2012synchronization}. In this work, we consider three popular models of coupled oscillators to study their synchronization behavior, described below.
	
	\subsection{Kuramoto model} The Kuramoto model \cite{strogatz2000kuramoto,acebron2005kuramoto} of coupled oscillators is perhaps one of the most well-studied models in the dynamical system community. Consider a graph $G = (V, E)$ and a continuous phase space $\Omega = \mathbb{R}/2\pi\mathbb{Z}$. The evolution of the phase dynamics of the initial phase configuration $X_0: V \rightarrow \Omega$ governed by the Kuramoto model of coupled oscillators is determined by the following system of ordinary differential equations in \eqref{eq:kuramoto}
	\begin{align}\label{eq:kuramoto}
		\frac{d}{dt} X_t(v) = \omega_v + K \sum_{u \in \mathcal{N}(v)} \sin(X_t(u) - X_t(v)) \quad\forall ~v \in V,
	\end{align}
	where $\mathcal{N}$(v) represents the set of neighboring nodes of $v$ in $G$, $\omega_v$ denotes the \textit{intrinsic frequency} of node $v$, and $K$ denotes the \textit{coupling strength} of the model. We discretize the ordinary differential equation in \eqref{eq:kuramoto} so that each `step' in Kuramoto dynamics is given by the following difference equation: 
	\begin{align}\label{eq:kura_discretize}
		X_{t+h}(v)-X_t(v) = h \left(\sum_{u\in \mathcal{N}(v)} K \sin(X_{t}(u)-X_{t}(v)) \right),
	\end{align}
	where we choose a step size of $h = 0.05$ and $K=1$.
	
\begin{figure}[ht!]
		\centering
		\includegraphics[width=1 \linewidth]{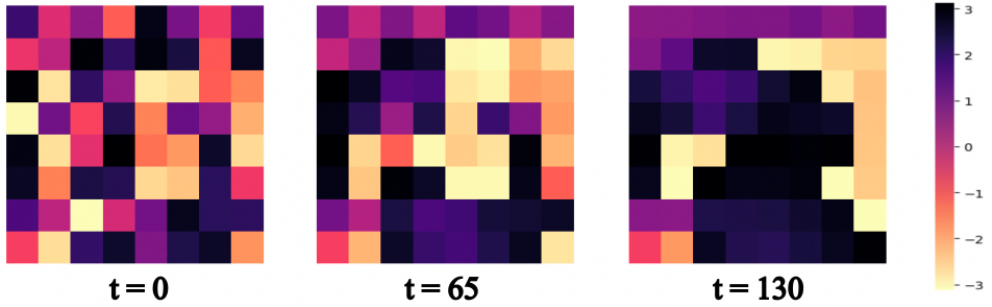}
		\caption{\textbf{Evolution of Kuramoto dynamics.} Representation of the evolution of the Kuramoto dynamics on an $8 \times 8$ 2D-grid graph with snapshots at time iterations $t=0, ~t=65, ~$ and $t=130$.}
		\label{fig:kura_example}
	\end{figure}

	In this paper, we assume the intrinsic frequency of all the nodes in $G$ are identical, or equivalently zero, by using a rotating frame of dynamics without loss of generality. We further assume the coupling strength to be unity. Note that synchronization is an absorbing state, in the sense that if $X_{\tau}$ is constant (e.g., synchronized), then $X_{t}$ is constant for all $t\ge \tau$. See Figure \ref{fig:kura_example} for an example of the Kuramoto dynamics evolving on an $8 \times 8$ grid.
	
	\subsection{The Firefly Cellular Automata}
	
\begin{figure}[ht!]
		\centering
		\includegraphics[width=1 \linewidth]{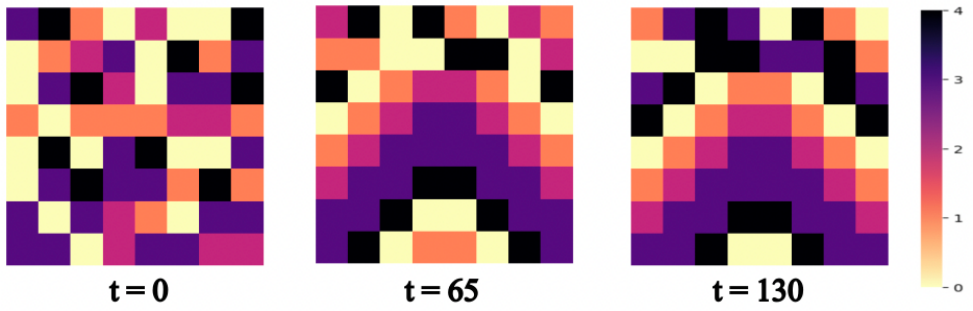}
		\caption{\textbf{Evolution of FCA dynamics.} Representation of the evolution of the FCA dynamics on an $8 \times 8$ 2D-grid graph with snapshots at time iterations $t=0, ~t=65, ~$ and $t=130$.}
		\label{fig:fca_example}
	\end{figure}

	The Firefly Cellular Automata (FCA) \cite{lyu2015synchronization,lyu2016synchronization} is a model to study discrete pulse-coupled oscillators. Consider a graph $G = (V, E)$ and $\kappa \ge 3$ to define $\Omega = \mathbb{Z}/\kappa\mathbb{Z}$ with an ordering $0 < 1 < \ldots < \kappa - 1$. The evolution of the phase dynamics of the initial phase configuration is governed by $X_0: V \rightarrow \Omega$ and we further define  $b(\kappa) = \lfloor\frac{\kappa-1}{2}\rfloor$ to be the \textit{blinking color} of the configuration. The time evolution of this $\kappa$-colored FCA dynamics is dictated by the update mapping $X \rightarrow X'$ in \eqref{eq:fca}
	\begin{align}\label{eq:fca}
		X'(v) = \begin{cases}
			X(v) & \parbox[t]{5.5cm}{\textup{if $X(v) > b(\kappa)$ and $v \in \mathcal{N}(v')$ such that $v' \in V$ has blinking state $b(\kappa)$}}\\
			X(v) + 1 & \textup{otherwise}
		\end{cases}
	\end{align}
	For all experiments in this work, we use the FCA model with $\kappa = 5.$ See Figure \ref{fig:fca_example} for an example of $5$-color FCA dynamics evolving on an $8 \times 8$ grid.
	
	\subsection{The Greenberg-Hastings model}
	
\begin{figure}[ht!]
		\centering
		\includegraphics[width=1 \linewidth]{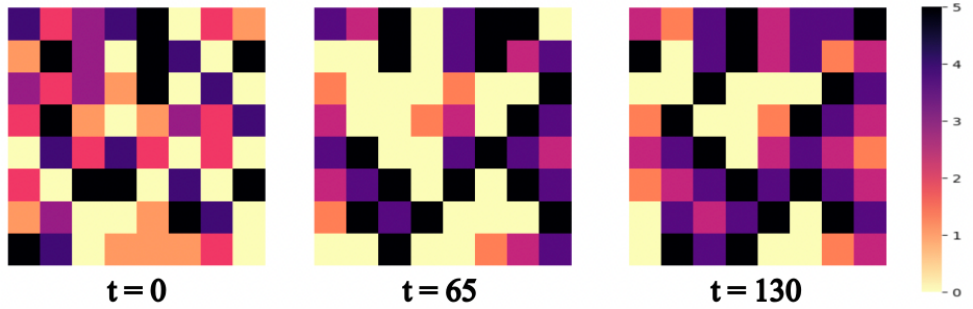}
		\caption{\textbf{Evolution of GHM dynamics.} Representation of the evolution of the GHM dynamics on an $8 \times 8$ 2D-grid graph with snapshots at time iterations $t=0, ~t=65, ~$ and $t=130$.}
		\label{fig:ghm_example}
	\end{figure}

	The Greenberg-Hastings model (GHM) \cite{greenberg1978spatial} is a popular model for studying discrete patterns of diffusion in excitable media. Consider a graph $G = (V, E)$ and define $\Omega = \mathbb{Z}/\kappa\mathbb{Z}$ with an ordering $0 < 1 < \ldots < \kappa - 1$. The evolution of the phase dynamics of the initial phase configuration with GHM is governed by $X_0: V \rightarrow \Omega.$ The time evolution of this $\kappa$-colored GHM dynamics is dictated by the mapping $X \rightarrow X'$ in \eqref{eq:ghm}
	\begin{align}\label{eq:ghm}
		X'(v) = \begin{cases}
			0 & \parbox[t]{5.5cm}{\textup{if $X(v) =0$ and $X(v') \neq 1 ~~\forall ~v' \in \mathcal{N}(v)$}}\\
			1 & \parbox[t]{5cm}{\textup{if $X(v) =0$ and $\exists ~v' \in \mathcal{N}(v)$ s.t. $X(v') = 1$}}\\
			X(v) + 1 & \textup{otherwise}
		\end{cases}
	\end{align}
	For all experiments in this work, we use the GHM model with $\kappa = 6.$ See Figure \ref{fig:ghm_example} for an example of $6$-color GHM dynamics evolving on an $8 \times 8$ grid.
	
	\section{Training Data, Models, Hyperparameters, and Experimental Details} \label{app:training_details}
	
	\begin{description}
		\item[Baseline Predictor] The baseline predictor is based on the concentration principle defined in Section \ref{sec:sufficient_conditions}, which predicts eventual synchronization of the evolving dynamics if the phase configuration is concentrated at any time during the observed dynamics, and flips an independent fair coin for prediction otherwise. 
		
		\item[Details for the subgraph-level experiments in Section \ref{sec:res_local}.] The data set we use for the subgraph-level experiments in Section \ref{sec:res_local} is generated as follows. From one of the three networks described in Table \ref{tab:networks}, we sample $k$-node subgraphs $F_{1},\dots,F_{N}$ with $N=10,000$ using the subgraph sampling oracle in Section \ref{sec:path_sampling}, where $k \in \{10, 15, 20, 25, 30\}$. 
		On the $i$th subgraph $F_{i}$, we randomly initialize and run $T'$ iterations of Kuramoto, FCA, and GHM dynamics, where $T'=200$, 100, and 100, respectively.  The data set consists of pairs $(\mathcal{X}_{i},y_{i})$, $i=1,\dots, N$, where $y_{i}$ is the indicator that the system is synchronized at time $T'$ and $\mathcal{X}_{i}$ is the $k\times k \times T$ CAT that encodes the dynamics on $F_{i}$ for the first $T<T'$ iterations, where $T=100$, $50$, and $8$ iterations for the Kuramoto, FCA and GHM dynamics, respectively. These $10,000$ data points are then uniformly randomly divided into training and testing sets, consisting of 80$\%$ and 20$\%$ of the examples, respectively. 
		
		For \textup{LLDM-T} in Table \ref{tab:acc_table}, we distill the generated data based on three observations from the literature as described in \ref{sec:sufficient_conditions}, namely graph density and initial half-circle concentration. In particular, we sample a set of subgraphs from our parent graphs and then distill the top $10\%$ densest and top $10\%$ sparsest, based on their edge density. In addition to these, we also select certain configurations of dynamics-network pairs where the dynamics follow the half-circle concentration, such that the total number of data points is $10,000$ which we again split into train and test sets (80\% and 20\%). 
		
		For \textup{LLDM} and \textup{LLDM-T} in Table \ref{tab:acc_table}, we used block minimization-type iterative algorithms for supervised matrix factorization \cite{lee2023supervised} and nonnegative matrix factorization algorithms \cite{lyu2023block} for $250$ iterations, respectively. 
		
		\item[Details for the global-level experiments in Section \ref{sec:res_global}.] 
		
		For the global-level experiments, we take the parent graph $G$ to be one of the 500 instances of the 300-node graphs in the dataset \dataset{NWS}'  in Table \ref{tab:networks}. These graphs are generated by the \dataset{NWS} model with the circulant graph with $k=12$ nearest-neighbors and $p = 0.4$ probability of adding a new edge independently between each non-adjacent pair of nodes. On each $G$, we simulate the Kuramoto and the $5$-color FCA dynamics with random initial configuration for $T'=50$ and $100$ iterations of FCA and Kuramoto dynamics, respectively. This creates 500 pairs of dynamics on $G$ and indicator $y_{G}$ of whether the system on $G$ is globally synchronized at time $T'$. We choose the initial configurations randomly so that half of them globally synchronize at time $T'$ (so $y_{G}=1$) and the other half do not (so $y_{G}=0$). These pairs are split into 80$\%$ training set and 20$\%$ testing set. From each $G$, we sample a single trajectory of 50 iterations of the MCMC $k$-path ($k \in \{10, 20, 30\})$ motif sampling algorithm (see Section \ref{sec:path_sampling}), which gives a sequence of $k$-node subgraphs $F_{1},\dots,F_{50}$. We then restrict the dynamics on $G$ on the sampled subgraphs, thereby creating pairs $(\X_{1},y_{G}),\dots,(\X_{50}, y_{G})$ of $k\times k\times T$ CATs and (global) synchronization indicator, where $T$ is varied between 10 and 100. In this way, we create a total of 500*50 pairs of observed subgraph dynamics and synchronization indicators. The block-coordinate descent algorithm for SMF \cite{lee2022supervised} is run for $250$ iterations to jointly learn the dictionary $\mathcal{D}$ (for ranks $R=\in \{4,25,100\}$) of latent dynamic filters and vector of regression coefficients $\Beta$ from the training data set. 
		
		\item[FFNN architecture.] The FFNN architecture we use for the experiments is one with four fully connected layers, where each intermediate layer has $100$ hidden nodes, batch normalization, and uses the ReLU \cite{agarap2018deep} activation function. Further, we use a dropout \cite{srivastava2014dropout} with $p=0.25$ on each of our layers to prevent model overfitting.

		\item[Supervision tuning parameter.] For our Supervised Matrix Factorization (SMF) based experiments (See Section \ref{sec:results} and Appendix \ref{app: SMF}), we report the best results on doing a grid search on hyperparameter $\xi$ (see \eqref{eq:SMF_1}) with choices $\xi \in [0.1, 0.5, 1.0].$ A higher value of $\xi$ indicates that the model will be penalized more for learning filters that do not reconstruct the original data well, but not so much for wrong label predictions, and vice-versa.
	\end{description}
	
	\section{Background on Feature extraction by NMF} \label{app: nmf}
	
\begin{figure}[ht!]
		\centering
		\includegraphics[width=1 \linewidth]{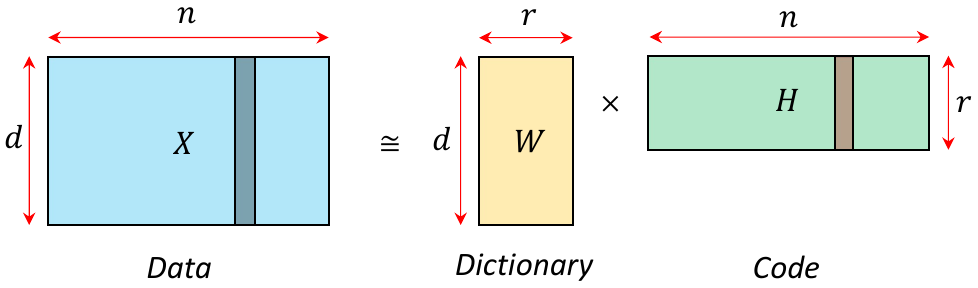}
		\caption{\textbf{Illustration of  matrix factorization.} Each column of the data matrix is approximated by the linear combination of the columns of the dictionary matrix with coefficients given by the corresponding column of the code matrix.}
		\label{fig:NMF_diagram}
	\end{figure}

	Suppose we are given with $n$ labeled signals $(\x_{i}, y_{i})$ for $i=1,\dots, n$, where $\x_{i}\in \R^{p}$ is a $p$-dimensional signal and $y_{i} \in \{ 0,1 \}$ is its binary label. Suppose  $p$ is large (e.g., high-dimensional signals) and there is a small number $R$ of latent feature vectors $\w_{1},\dots,\w_{R}$, forming a `dictionary' matrix $\W\in \R^{p\times R}$ such that each high-dimensional signal $\x_{i}$ can be approximated by some linear combination $\W \h_{i}$ for some $\h_{i}\in \R^{R}$. In this way, the high-dimensional signal $\x_{i}$ is compressed into a low-dimensional signal $\h_{i}$. The problem of finding suitable factor matrices $\W$ and $\mathbf{H}=[\h_{1},\dots,\h_{R}]$ can be formulated by a \textit{matrix factorization} problem $\mathbf{X}\approx \W \mathbf{H}$. That is, each column of the data matrix is approximated by the linear combination of the columns of the \textit{dictionary matrix} $W$ with coefficients given by the corresponding column of the \textit{code matrix} $H$ (see Figure \ref{fig:NMF_diagram}). Variants of matrix factorization problems have been investigated under many names over the decades, each with different assumptions and constraints: dictionary learning, factor analysis, topic modeling, and component analysis. It has applications in text analysis, image reconstruction, medical imaging, bioinformatics, network dictionary learning, and many other scientific fields more generally \cite{sitek2002correction, berry2005email, berry2007algorithms, chen2011phoenix, taslaman2012framework, boutchko2015clustering, ren2018non, lyu2021learning}.
	
	Along with Principal Component Analysis (PCA) \cite{wold1987principal}, nonnegative Matrix Factorization (NMF)\cite{lee1999learning} is a classical matrix factorization setting where both factor matrices $W$ and $H$ are constrained to be nonnegative. In the simplest form, NMF is formulated as the following bi-convex, constrained optimization problem. Given a data matrix $\mathbf{X}\in \mathbb{R}^{d\times n}$, the goal is to find a nonnegative dictionary $\mathbf{W}\in \mathbb{R}^{d\times r}_{\ge 0}$ and nonnegative code matrix $\mathbf{H}\in \mathbb{R}^{r\times n}_{\ge 0}$ by solving the following optimization problem:
	\begin{align}\label{eq:NMF_error1}
		\inf_{\substack{\W\in \mathbb{R}^{d\times r}_{\ge 0},\,   \mathbf{H}\in \mathbb{R}^{r\times n}_{\ge 0} }}  \lVert \mathbf{X} - \W\mathbf{H}  \rVert_{F}^{2}, 
	\end{align}
	where $\lVert A \rVert_{F} = \sqrt{\sum_{i,j} A_{ij}^{2}}$ denotes the matrix Frobenius norm.
	
	A consequence of the nonnegativity constraint on the code matrix is that one must represent the columns of the data matrix only using nonnegative linear combinations of the columns of the dictionary matrix $\W$. Since the columns of $\W$ are also constrained to be nonnegative, every feature captured in columns of $\W$ can only additively (rather than subtractively) be combined to explain the data points (columns in $\mathbf{X}$). This allows one to interpret the columns of $\W$ as `parts' of the data and the columns of $\mathbf{H}$ as their `contribution' in composing the columns of the data. 
	
	\section{Background on Supervised Matrix Factorization (SMF)}\label{app: SMF}
	
	Note that the NMF formulation in \eqref{eq:NMF_error1} does not incorporate the labels $y_{1},\dots,y_{n}$. This means the dictionary matrix $\W$ is for the best possible reconstruction of the data matrix $\mathbf{X}$, but the best reconstructive dictionary $\W$ may not be very effective for the classification tasks. 
	
	In the \textit{supervised matrix factorization} (SMF) literature \cite{mairal2008supervised, lee2023supervised}, one desires a dictionary that is reconstructive as well as \textit{discriminative} in that such a compressed representation of signals is adapted to predicting the class labels $y_{i}$. In order to learn a dictionary matrix $\W$ that is both data-reconstructive and label-discriminative, we jointly model the pairs $(\x_{i},y_{i})$ of high-dimensional signal $\x_{i}$ and binary label $y_{i}$ as  
	\begin{align}\label{eq:SMF_model}
		\x_{i} \approx \W\h_{i} \quad  \text{and} \quad y_{i}\,|\, \x_{i}\sim \textup{Bernoulli}\left( p_{i} \right), 
	\end{align}
	where $\W\in\R^{p\times R}_{\ge 0}$ is some unknown \textit{nonnegative} dictionary matrix and the $p_{i}$ is the predictive probability given by 
	\begin{align}\label{eq:SMF_p_def}
		p_{i}:= \frac{\exp(\langle \Beta, \W^{T}\x_{i} \rangle)}{1+\exp(\langle \Beta, \W^{T}\x_{i} \rangle}
	\end{align}
	for $i=1,\dots,n$. 
	Here we can interpret the product $\W^{T}\x_{i}$ as performing a convolution on the $p$-dimensional signal $\x_{i}$ using the $R$ columns in $\W$. Since $\W$ is nonnegative, such convolution computes the proximity score of the pattern in each column of $\W$ observed in observed signal $\x_{i}$. In this sense, we view the columns of $\W$ as `filters' that encode $R$ particular patterns we seek to detect in $\x_{i}$. The vector $\W^{T}\x_{i}\in \R^{R}$ of proximity scores is then used as input to the logistic classifier with regression coefficients in $\Beta$. If the $j$th coordinate of $\Beta$ is positive (resp., negative), then the proximity score $\langle \w_{j}, \x_{i}\rangle$ of the signal $\x_{j}$ with $j$th filter $\w_{j}$ is positively (resp., negatively) associated with $y_{i}$ being one (resp., zero). 
	
	We formulate the following joint optimization problem for fitting \textit{nonnegative} SMF model to the data-label pairs: 
	\begin{align}\label{eq:SMF_1}
		&\min_{\W\ge 0, \mathbf{H}\ge 0, \Beta} \,\,  \left( \sum_{i=1}^{n} \ell(y_{i}, \langle \Beta,\,  \W^{T} \x_{i} \rangle)  +   \xi \lVert \mathbf{X} - \W\mathbf{H}\rVert_{F}^{2} \right),
	\end{align}
	where $\W\in \R^{p\times R}_{\ge 0}$, $\mathbf{H}\in \R^{R\times n}_{\ge 0}$, $\Beta\in \R^{R}$, and the negative log-likelihood function $\ell$ is defined by 
	\begin{align}\label{eq:ell_log_likelihood}
		\ell(y_{i}, \langle \Beta,\,  \W^{T} \x_{i} \rangle)
		&:= - y_{i}\log p_{i}  - (1-y_{i}) \log (1-p_{i}) 
	\end{align}
	with $p_{i}$ as in \eqref{eq:SMF_p_def}. Here, the \textit{tuning parameter} $\xi\ge 0$ controls the trade-off between the two objectives of classification and matrix factorization.
	
	The objective function in \eqref{eq:SMF_1}, say $F(\W,\mathbf{H}, \Beta)$, is convex in each of the three variables $\W,\mathbf{H}$, and $\Beta$ while the other two variables are held fixed. Hence we can employ cyclic block coordinate descent (BCD) algorithms\cite{lyu2023block}. In particular, the simplest cyclic block minimization algorithm for solving \eqref{eq:SMF_1} reads as 
	\begin{align}\label{eq:SBM_BCD}
		\begin{cases}
			\W_{k+1} &\leftarrow    \argmin_{\W\ge 0} \,\, F(\W,\mathbf{H}_{k},\Beta_{k}) \\
			\mathbf{H}_{k+1} &\leftarrow    \argmin_{\mathbf{H}\ge 0} \,\, F(\W_{k+1},\mathbf{H},\Beta_{k}) \\
			\Beta_{k+1} &\leftarrow    \argmin_{\Beta} \,\, F(\W_{k+1},\mathbf{H}_{k+1},\Beta).
		\end{cases}
	\end{align}
	Each sub-problem in \eqref{eq:SBM_BCD} is a convex optimization problem, which can be solved by using standard algorithms such as projected gradient descent\cite{boyd2004convex}. See Lee, Lyu, and Yao \cite{lee2022supervised} for a more detailed implementation of BCD for nonnegative SMF \eqref{eq:SMF_1} and convergence guarantees.

	\section{Performance Comparison for Various Graph Embedding Methods}\label{app:graph_embedding}
	
	In this section, we provide additional experiments in order to demonstrate that the way we encode the topological features of the underlying graph does not make a significant difference in the prediction accuracy for the synchronization prediction of coupled oscillators on these graphs. 
	
	Graph embedding techniques like \texttt{DeepWalk} \cite{perozzi2014deepwalk}, \texttt{graph2vec} \cite{narayanan2017graph2vec}, and Spectral Embedding \cite{ng2002spectral} have been used extensively in recent years to embed a graph into a low-dimensional vector space while preserving the structure of the graph. In this paper, we have used a rather simple colored adjacency tensor (CAT) \cite{eq:CAT} to encode $T$-iterations of dynamics on a $k$-node subgraph into a $k\times k\times T$ nonnegative tensor. As the name suggests, it is a stack of the adjacency matrix of the underlying graph, with additional information on the phase difference along the edges at each time. One may wonder if one uses a more sophisticated graph embedding algorithm to encode a graph-dynamics pair, then one would get a potentially significant performance gain in synchronization prediction problems. However, we argue that this is not the case. 
	
	\texttt{Spectral embedding} is a classical graph embedding technique that uses top eigenvectors of the graph Laplacian matrix as low-dimensional vector representations of the nodes of a network. The objective of the \texttt{DeepWalk} is to learn a mapping of nodes into a low dimensional Euclidean space such that two nodes that co-appear frequently in random walk sequences on the network would be assigned with their vector representations with the large inner product; two nodes that do not co-appear frequently will be nearly orthogonal after the embedding. The objective of the \texttt{graph2vec} is to learn low dimensional graph embeddings in an unsupervised manner, primarily used for graph classification. 
	
	In Table \ref{table:embed}, we perform synchronization prediction on $k=20$-node subgraphs from \dataset{UCLA} network using a logistic classifier on data sets generated by using various different methods to encode the graph topology. The data generation setting is identical to that for the subgraph-level experiments in Section \ref{sec:res_local}, which we explained in detail in Appendix \ref{app:training_details}. In this table, `dynamics' means the $k\times T$ matrix encoding of the $T$-iterations of dynamics on the graph. We append this matrix with four different encodings of the underlying graph: \texttt{Adjacency matrix}, \texttt{spectral embedding}, \texttt{DeepWalk}, and \texttt{graph2vec}. 
	
	\begin{table}[!ht]
		\label{table:embed}
		\caption{Effect on Logistic Regression Prediction Accuracy of Adding Graph Embedding Features in Addition to Dynamics.}
		\begin{ruledtabular}
			\begin{tabular}{c|ccc}
				\textbf{Embedding Technique} & Kuramoto & FCA & GHM\\
				\hline
				\textup{Dynamics} & 96.4\% & 92.7\% & 90.6\%\\
				\textup{Dynamics +} \texttt{Adjacency Matrix} & 96.3\% & 92.9\% & 91.4\%\\
				\textup{Dynamics +} \texttt{Spectral Embedding} & \textbf{96.9}\% & 92.7\% & 90.9\%\\
				\textup{Dynamics +} \texttt{DeepWalk} & 96.5\% & 92.1\% & 91.3\% \\
				\textup{Dynamics +} \texttt{graph2vec} & 96.8\% & \textbf{93.1}\% & \textbf{91.8}\% 
			\end{tabular}
		\end{ruledtabular}
	\end{table}
	
	We observe that encoding the underlying graph using various methods does not seem to provide the model with much additional information. The highest accuracy gain even for modern embedding methods compared to providing the model with only dynamics and the adjacency matrix is in the case of Kuramoto, with a $0.6\%$ gain, which does not lead to a significant difference. Moreover, this trend remains relevant for all three coupled oscillator models. Overall, changing the graph embedding technique seems to have little to no effect on model performance, which is why a simple adjacency-matrix-based canonical representation that is encoded by our CATs (See Figure \ref{fig:CAT}) is already sufficient to provide the model with enough graph topology information.
	
	\vspace{0.5cm}
	
	\section{Extended subgraph-level prediction accuracies}\label{app:extended_subgraph_results}
	
	In this section, we provide the full set of results of prediction accuracy of various methods versus LLDM on the subgraph level for $k$-node subgraphs where $k \in \{10, 15, 20, 25, 30\}$ with the FCA, Kuramoto, and GHM dynamics induced on them. Table \ref{tab:acc_table_nws} represents accuracies for subgraphs sampled from NWS, table \ref{tab:acc_table_caltech} represents accuracies for subgraphs sampled from \dataset{Caltech}, and table \ref{tab:acc_table_ucla} represents accuracies for subgraphs sampled from \dataset{UCLA}. Each accuracy metric reported is the mean accuracy from $5$ seeds, based on a grid search across tuning parameter $\xi \in [0.1, 0.5, 1.0].$
	
	\begin{table*}[!ht]
		\label{tab:acc_table_nws}
		\caption{
			Prediction accuracy of various methods for FCA, Kuramoto, and GHM dynamics on subgraphs with $k$ number of nodes where $k \in \{10, 15, 20, 25, 30\}$ sampled from a large-connected \dataset{NWS} parent graph. All accuracy values are an average of $5$ seeds. The highest accuracy for each setting is indicated in \textbf{bold}. Whenever the average values of accuracy are equal for two methods, \textbf{both} are represented with \textbf{bold} font.}
		\begin{ruledtabular}
			\begin{tabular}{c|ccccc|ccccc|cccccc}
				&&&FCA&&&
				&&Kuramoto&&&
				&&GHM&&&\\
				\hline
				&$k=10$&$k=15$&$k=20$&$k=25$&$k=30$
				&$k=10$&$k=15$&$k=20$&$k=25$&$k=30$
				&$k=10$&$k=15$&$k=20$&$k=25$&$k=30$\\
				\hline
				Baseline&80.2&78.2&69.6&71.0&69.4&56.5&63.2&66.1&69.9&64.3&89.6&78.7&76.2&69.6&65.6\\
				LogReg&92.7&\textbf{95.2}&94.8&89.1&76.4&83.5&81.8&80.1&82.4&80.2&95.5&92.7&\textbf{93.6}&90.1&92.8\\
				FFNN&92.2&94.9&\textbf{95.4}&\textbf{90.2}&82.8&\textbf{84.2}&\textbf{82.5}&\textbf{83.4}&83.8&79.2&94.4&95.8&91.8&\textbf{91.4}&\textbf{94.1}\\
				\hline
				LLDM-T (R=2)&82.7&89.4&81.3&85.0&75.9&74.3&76.6&76.0&77.1&80.1&85.5&90.2&89.0&88.0&89.7\\
				LLDM-T (R=8)&86.5&91.3&81.9&86.0&75.7&75.2&77.5&77.1&77.4&\textbf{80.3}&81.6&90.4&81.7&78.1&83.0\\
				LLDM (R=2)&91.5&92.2&93.2&87.0&77.2&81.4&76.1&77.5&77.2&78.3&95.9&93.3&91.0&89.4&89.5\\
				LLDM (R=8)&\textbf{93.0}&92.8&94.2&88.1&\textbf{84.8}&82.1&81.6&78.4&\textbf{84.4}&78.4&\textbf{96.0}&\textbf{96.2}&92.3&89.0&\textbf{94.1}\\
			\end{tabular}
		\end{ruledtabular}
	\end{table*}
	
	\begin{table*}[!ht]
		\caption{\label{tab:acc_table_caltech}
			Prediction accuracy of various methods for FCA, Kuramoto, and GHM dynamics on subgraphs with $k$ number of nodes where $k \in \{10, 15, 20, 25, 30\}$ sampled from a large-connected \dataset{Caltech} parent graph. All accuracy values are an average of $5$ seeds. The highest accuracy for each setting is indicated in \textbf{bold}. Whenever the average values of accuracy are equal for two methods, \textbf{both} are represented with \textbf{bold} font.}
		\begin{ruledtabular}
			\begin{tabular}{c|ccccc|ccccc|cccccc}
				&&&FCA&&&
				&&Kuramoto&&&
				&&GHM&&&\\
				\hline
				&$k=10$&$k=15$&$k=20$&$k=25$&$k=30$
				&$k=10$&$k=15$&$k=20$&$k=25$&$k=30$
				&$k=10$&$k=15$&$k=20$&$k=25$&$k=30$\\
				\hline
				Baseline&78.4&74.3&68.8&68.3&63.7&52.5&61.6&68.3&66.1&65.4&85.3&74.5&72.6&66.7&65.8\\
				LogReg&91.8&92.2&90.5&91.5&\textbf{89.4}&82.6&80.2&81.2&85.6&81.7&\textbf{96.5}&93.2&92.9&88.3&91.5\\
				FFNN&\textbf{92.9}&\textbf{93.8}&\textbf{91.4}&\textbf{92.2}&87.2&83.6&\textbf{81.2}&80.8&\textbf{90.4}&\textbf{81.9}&95.3&\textbf{95.2}&\textbf{95.7}&90.8&\textbf{91.8}\\
				\hline
				LLDM-T (R=2)&89.0&88.8&78.4&81.9&76.2&\textbf{83.9}&78.2&79.5&82.4&81.0&96.0&92.5&89.2&86.4&84.3\\
				LLDM-T (R=8)&86.9&90.6&80.1&80.3&82.4&82.7&78.6&79.8&84.8&81.7&92.9&92.4&84.9&70.9&76.7\\
				LLDM (R=2)&88.8&91.4&81.1&84.6&80.4&82.6&76.8&79.2&83.7&80.1&95.4&92.2&89.0&84.7&82.4\\
				LLDM (R=8)&90.4&93.0&89.5&87.7&87.8&82.8&77.9&\textbf{82.2}&89.0&79.6&95.5&93.4&93.8&\textbf{91.6}&90.1\\
			\end{tabular}
		\end{ruledtabular}
	\end{table*}
	
	\begin{table*}[!ht]
		\caption{\label{tab:acc_table_ucla}
			Prediction accuracy of various methods for FCA, Kuramoto, and GHM dynamics on subgraphs with $k$ number of nodes where $k \in \{10, 15, 20, 25, 30\}$ sampled from a large-connected \dataset{UCLA} parent graph. All accuracy values are an average of $5$ seeds. The highest accuracy for each setting is indicated in \textbf{bold}. Whenever the average values of accuracy are equal for two methods, \textbf{both} are represented with \textbf{bold} font.}
		\begin{ruledtabular}
			\begin{tabular}{c|ccccc|ccccc|cccccc}
				&&&FCA&&&
				&&Kuramoto&&&
				&&GHM&&&\\
				\hline
				&$k=10$&$k=15$&$k=20$&$k=25$&$k=30$
				&$k=10$&$k=15$&$k=20$&$k=25$&$k=30$
				&$k=10$&$k=15$&$k=20$&$k=25$&$k=30$\\
				\hline
				Baseline&80.6&77.4&70.2&71.2&68.3&56.4&65.3&69.2&68.4&64.9&86.5&76.3&71.4&67.2&65.2\\
				LogReg&89.8&89.2&94.1&\textbf{93.8}&91.2&91.3&90.2&91.4&92.8&94.6&\textbf{94.2}&\textbf{92.4}&\textbf{95.6}&92.3&92.3\\
				FFNN&90.4&\textbf{91.1}&\textbf{95.6}&93.2&\textbf{92.5}&\textbf{92.8}&\textbf{91.3}&\textbf{92.4}&93.6&\textbf{95.2}&93.8&92.2&94.8&\textbf{93.4}&92.9\\
				\hline
				LLDM-T (R=2)&77.9&79.7&80.4&82.3&80.1&79.8&84.1&87.1&83.2&87.2&92.3&92.2&89.4&92.5&91.4\\
				LLDM-T (R=8)&78.3&80.9&84.7&87.3&84.1&82.6&85.8&87.6&86.5&90.2&94.6&91.1&87.9&89.5&88.4\\
				LLDM (R=2)&86.8&86.5&90.5&90.8&86.4&79.1&86.2&85.9&85.2&93.5&93.2&91.6&90.7&91.1&91.2\\
				LLDM (R=8)&\textbf{91.1}&88.6&92.6&92.7&90.0&89.3&89.3&90.0&\textbf{94.2}&94.7&94.6&91.6&93.3&93.1&\textbf{95.9}\\
			\end{tabular}
		\end{ruledtabular}
	\end{table*}
	
	\section{latent dynamics filters or Dictionary Plots for Various Settings}\label{appendix:dictionary-plots}
	
\begin{figure*}[!ht]
		\centering
		\includegraphics[width=0.75\linewidth]{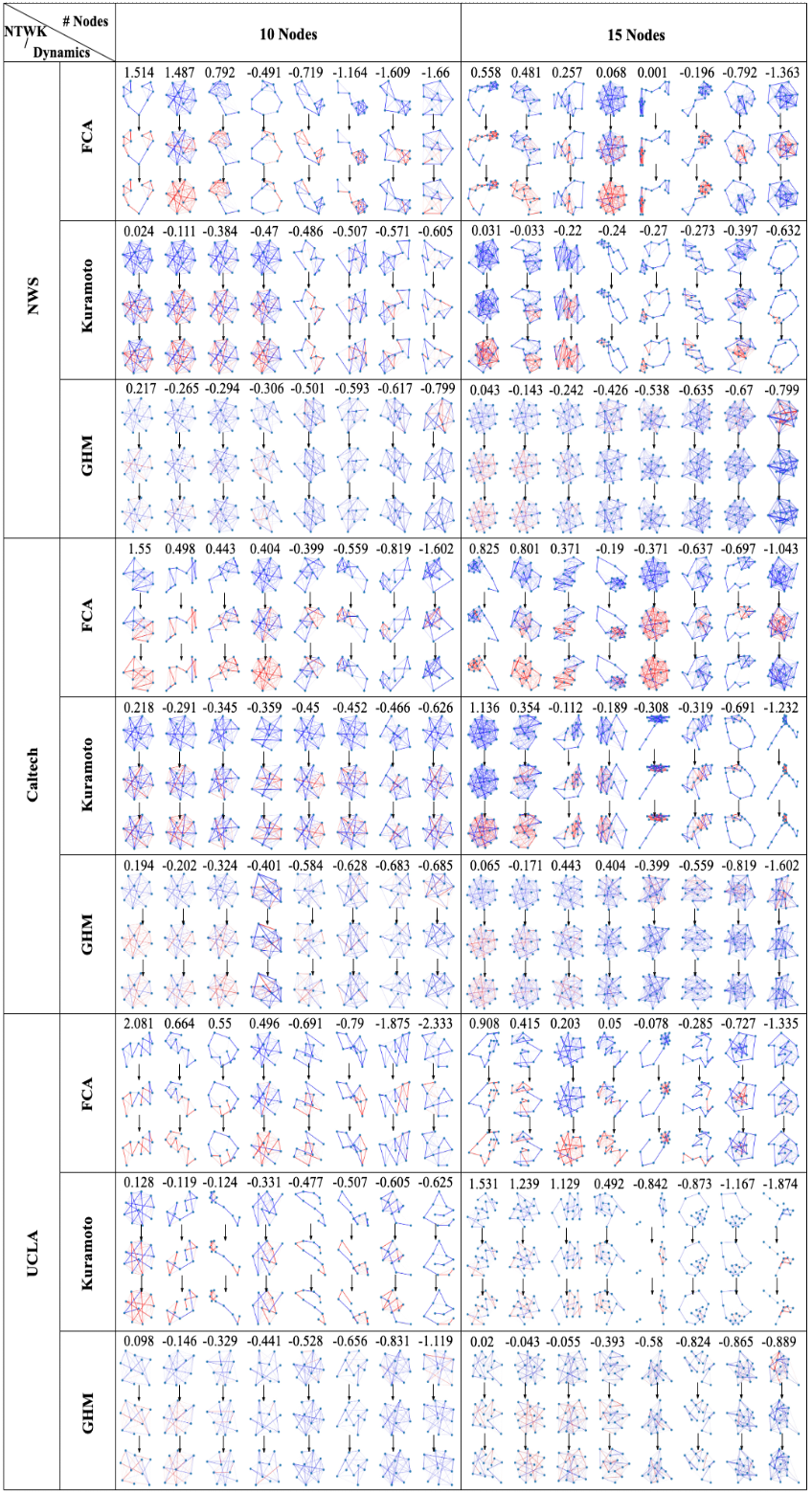}\caption{The 8-element SMF latent dynamics filters of FCA and Kuramoto on \dataset{Caltech}, \dataset{UCLA}, \dataset{NWS} networks of 10 nodes and 15 nodes}
		\label{fig:filters_appendix1}
	\end{figure*}

\begin{figure*}[!ht]
		\centering
		\includegraphics[width=0.80\linewidth]{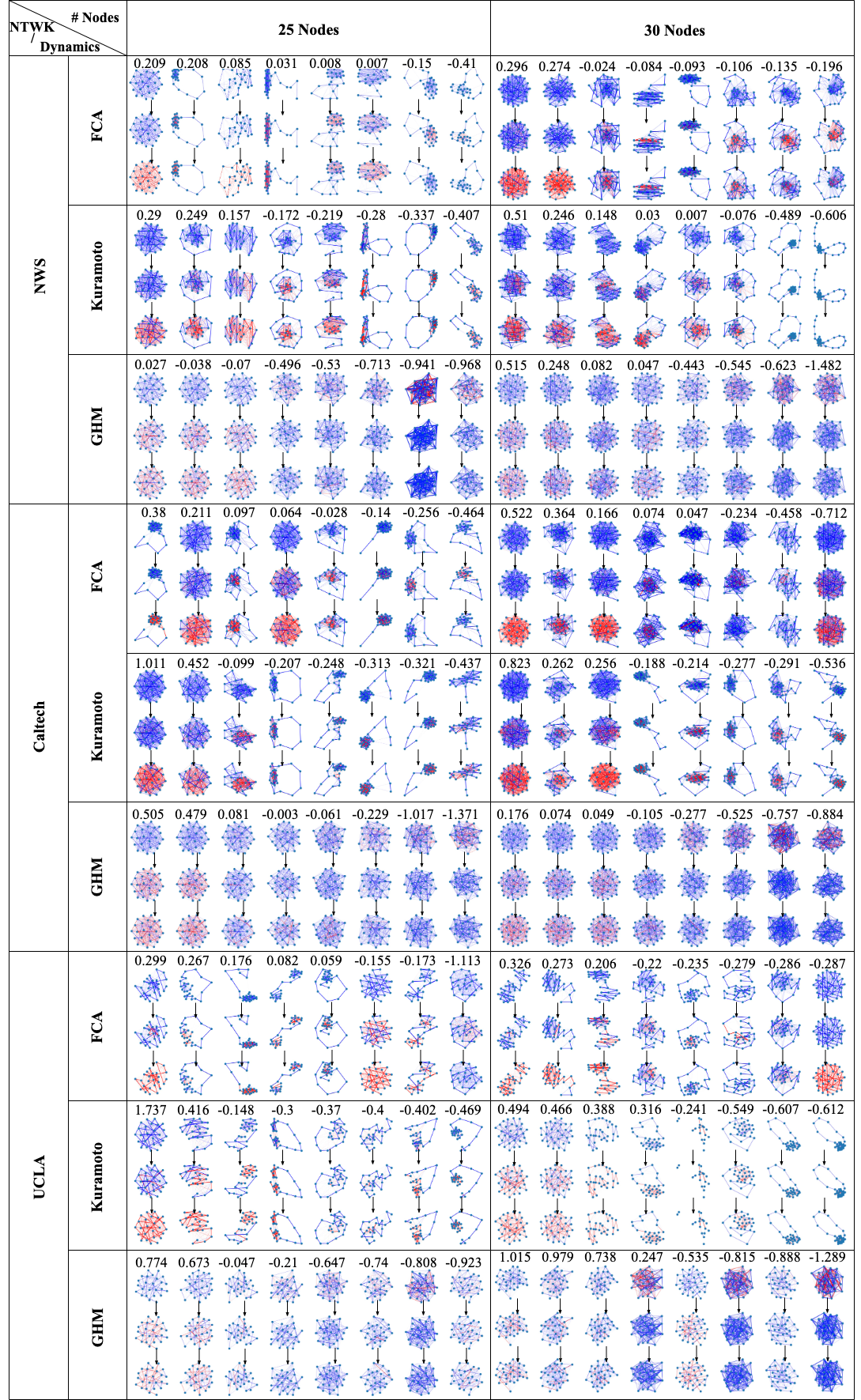}\caption{The 8-element SMF latent dynamics filters of FCA and Kuramoto on \dataset{Caltech}, \dataset{UCLA}, \dataset{NWS} networks of 25 and 30 nodes}
		\label{fig:filters_appendix2}
	\end{figure*}

	In this section, we show the rank $8$ dictionary atoms learned by LLDM from the three parent networks and three dynamics networks for subgraph sizes $k \in \{10, 15, 25, 30\}$ with $k=20$ in the main text Section \ref{sec:res_interpret}. See Figures \ref{fig:filters_appendix1} and \ref{fig:filters_appendix2}. 
	
	\section{Goodness-of-fit Plots}\label{app:gof_extended}
	
\begin{figure*}[!ht]
		\centering
		\includegraphics[width=0.45\textwidth]{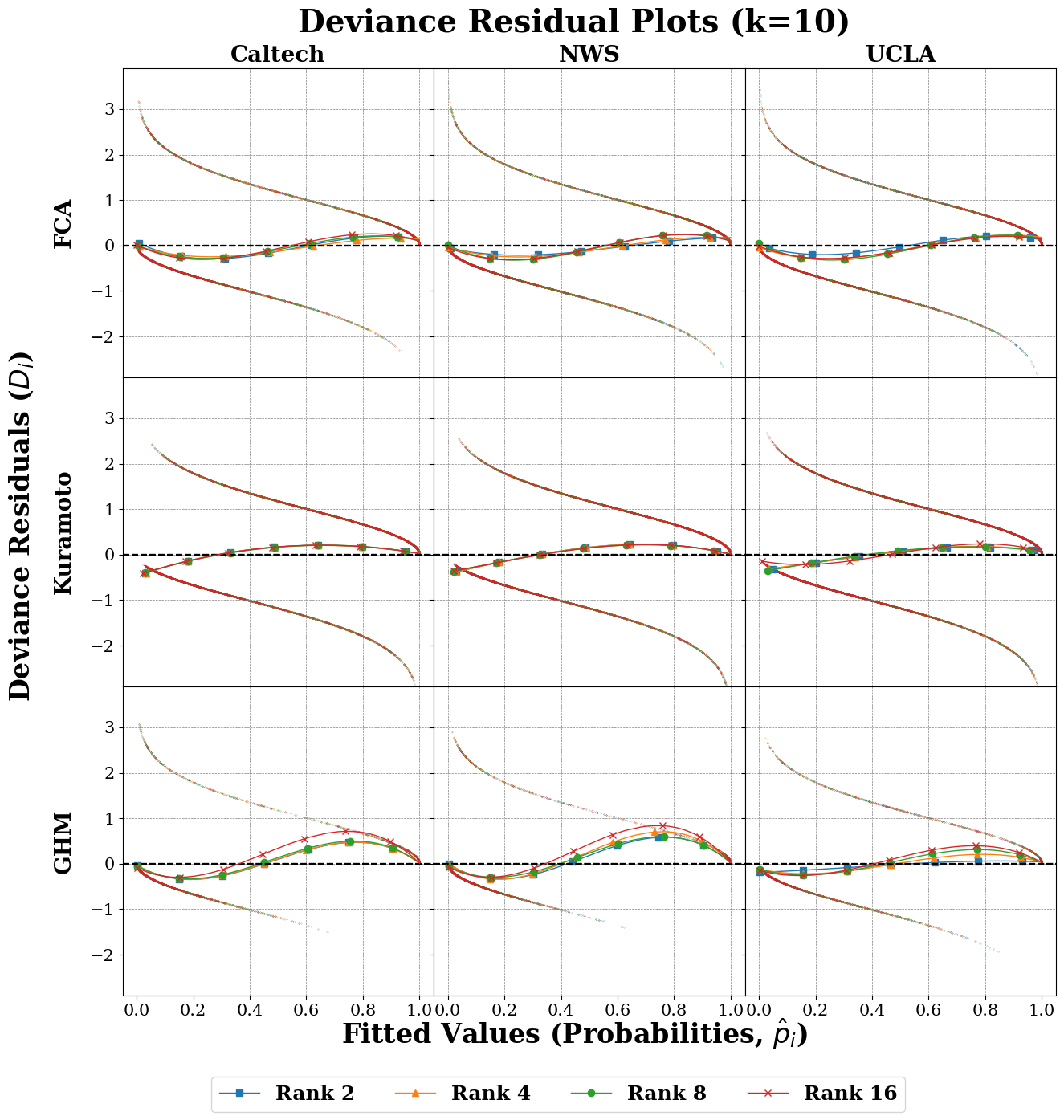}
		\includegraphics[width=0.45\textwidth]{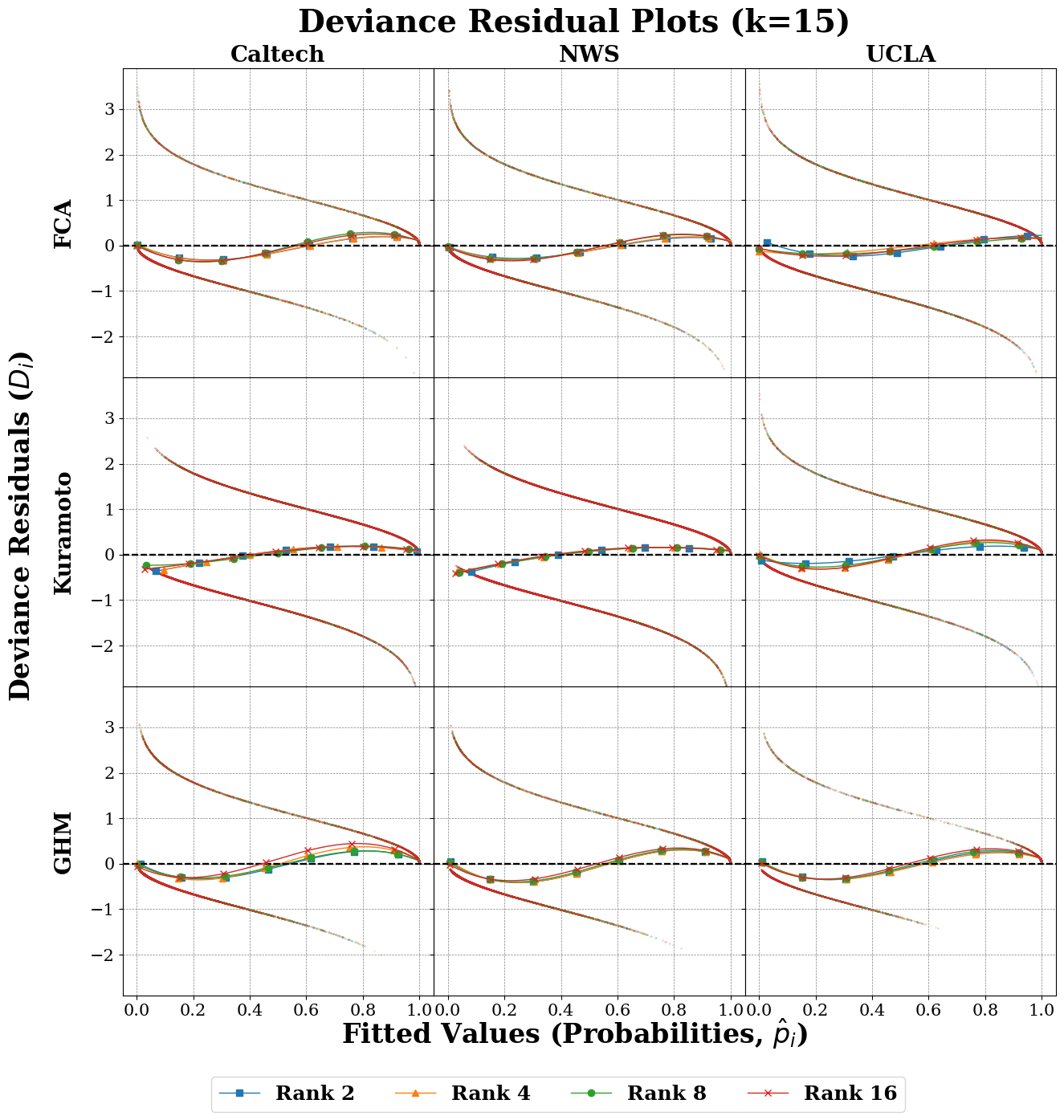}
		\includegraphics[width=0.45\textwidth]{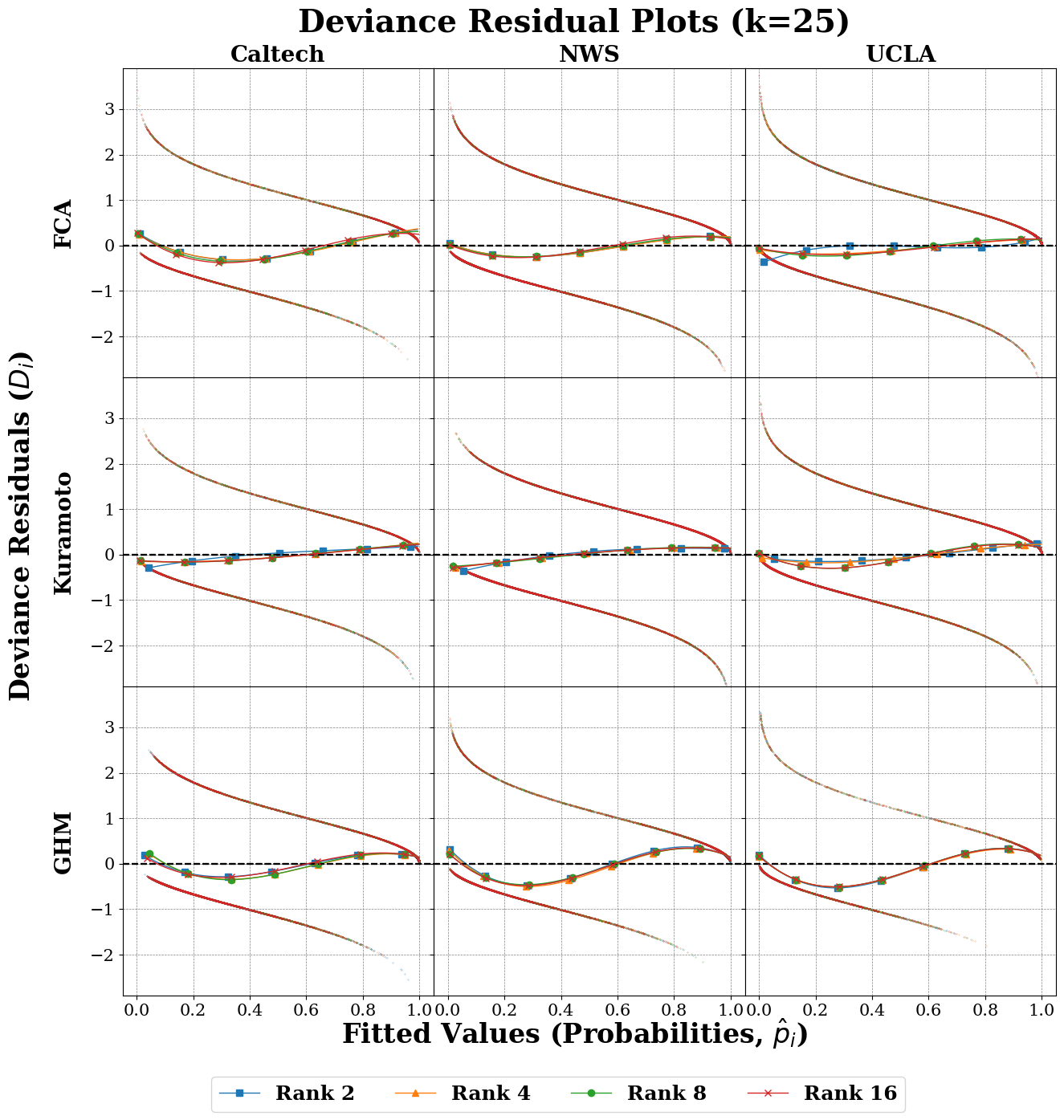}
		\includegraphics[width=0.45\textwidth]{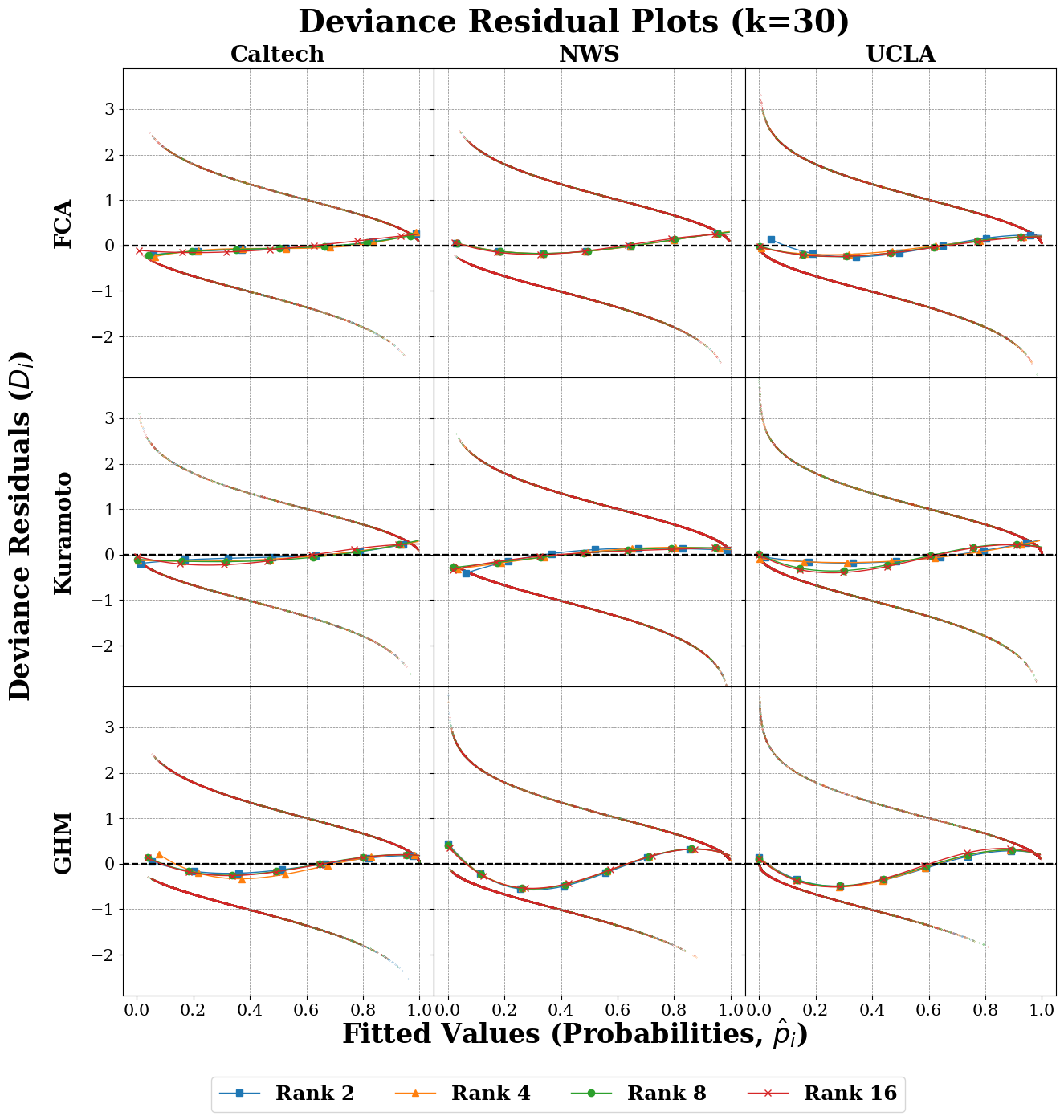}
		\caption{Goodness-of-fit with deviance residuals for $k=10,\;15,\;25$ and 30 on rank 2, 4, 8 and 16 for \dataset{Caltech}, \dataset{UCLA}, and \dataset{NWS} graphs with FCA, Kuramoto and GHM models}
		\label{fig:dev15}
	\end{figure*}

	In this section, we show the rank $2, 4, 8,$ and $16$ goodness-of-fit plots in Figure \ref{fig:dev15} with deviance residuals for subgraph sizes $k \in \{10, 15, 25, 30\}$ with $k=20$ in the main text Section \ref{sec:goodness-of-fit}.
	
\end{document}